\documentclass[10pt]{article}
\usepackage[latin1]{inputenc}
\usepackage{epsfig}
\usepackage{color}
\usepackage[british,english]{babel}
\usepackage{amsthm}
\usepackage{amsmath}
\usepackage{amsfonts}
\usepackage{amssymb}
\usepackage{graphicx}
\def \noame{\noalign{\medskip}}
\newtheorem{corollary}{Corollary}[section]
\newtheorem{definition}[corollary]{Definition}

\newtheorem{lemma}[corollary]{Lemma}
\newtheorem{proposition}[corollary]{Proposition}

\newtheorem{theorem}[corollary]{Theorem}
\newfont{\sBlackboard}{msbm10 scaled 900}

\newcommand{\mylabel}[1]{\label{#1}
            \ifx\undefined\stillediting
            \else \fbox{$#1$}\fi }
\newcommand{\BE}{\begin{equation}}

\newcommand{\EEQ}{\end{equation}}
\newcommand{\rfb}[1]{\mbox{\rm
   (\ref{#1})}\ifx\undefined\stillediting\else:\fbox{$#1$}\fi}

\newfont{\Blackboard}{msbm10 scaled 1200}

\newfont{\roma}{cmr10 scaled 1200}

\def\CC{\rm \hbox{C\kern-.56em\raise.4ex
         \hbox{$\scriptscriptstyle |$}\kern+0.5 em }}

\newcommand{\ep}{\varepsilon}
\newcommand{\rot}{{\rm rot}}

\def\n{|\kern -.05cm{|}\kern -.05cm{|}}

\newcommand{\mm}    {{\hbox{\hskip 0.5pt}}}

\newcommand{\bluff} {{\hbox{\raise 15pt \hbox{\mm}}}}
\usepackage{fancyhdr}

\makeatletter
\def\section{\@startsection {section}{1}{\z@}{-3.5ex plus -1ex minus
    -.2ex}{2.3ex plus .2ex}{\large\bf}}
\makeatother
\def\be{\begin{equation}}
\def\ee{\end{equation}}

\date{ }
\begin{document}
\thispagestyle{empty}
\title{\bf Homogenization of a micropolar fluid past a porous media with non-zero spin boundary condition}
\maketitle
%non-zero boundary conditions for microrotations}\maketitle
\vspace{-50pt}
%\author{ \center  Mar\'ia ANGUIANO\\
%Departamento de An\'alisis Matem\'atico. Facultad de Matem\'aticas.\\
%Universidad de Sevilla, P. O. Box 1160, 41080-Sevilla (Spain)\\
%anguiano@us.es\\}
\begin{center}
Francisco Javier SU\'AREZ-GRAU\footnote{Departamento de Ecuaciones Diferenciales y An\'alisis Num\'erico. Facultad de Matem\'aticas. Universidad de Sevilla. 41012-Sevilla (Spain) fjsgrau@us.es}
 \end{center}

\vskip20pt

 \renewcommand{\abstractname} {\bf Abstract}
\begin{abstract} 
We consider a micropolar fluid flow  in a media perforated by periodically distributed obstacles of size $\ep$.  A  non-homogeneous boundary condition for microrotation is considered: the microrotation is assumed to be proportional to
the rotation rate of the velocity on the boundary of the obstacles. The existence and uniqueness of solution is
analyzed. Moreover, passing to the limit when $\ep$ tends to zero, an analogue of the classical micropolar Darcy law  in the theory of porous media is derived. \end{abstract}
\bigskip\noindent
 {\small \bf AMS classification numbers:}  76A05, 76M50, 76S05, 35B27, 35Q35. 

\section {Introduction}\label{S1}Micropolar fluid model is a non-Newtonian model which represents a generalization of the well-established Newtonian Navier-Stokes model which takes into account the microstructure of the fluid. It describes the behavior of numerous real fluids (e.g. polymeric suspensions, liquid crystals, muddy fluids, animal blood, etc.) better than the classical one. The related mathematical model expresses the balance of momentum, mass and angular momentum.  Thus, a new unknown function $\hat w$ called microrotation  (i.e. the angular velocity field of rotation of particles) is added to the classical pressure $\hat p$ and velocity $\hat u$. Consequently, Newtonian Navier-Stokes equations become coupled with a new vector equation coming from the conservation of angular momentum, see Eringen \cite{Eringen} and Lukaszewicz \cite{Luka}. In view of its application in porous media,
we can assume a small Reynolds number and neglect the nonlinear terms and so the following micropolar equations are considered
\begin{eqnarray}
-(\nu+\nu_r)\Delta \hat u+\nabla \hat p =2\nu_r{\rm rot}(\hat w )+\hat  f,  &\nonumber\\
\noame
{\rm div}(\hat u)=0,&\nonumber\\
\noame
-(c_a+c_d)\Delta\hat  w +4\nu_r\hat w =2\nu_r{\rm rot}(\hat u) + \hat g. &\nonumber
\end{eqnarray}
The different viscosities  $\nu$, $\nu_r$, $c_a$ and $c_d$ are strictly positive and characterize the isotropic properties of the fluid,  $\hat f$  represents the body force and  $\hat g$  the body torque.

Solution to the governing equations  in the presence of solid boundaries requires imposing appropriate boundary conditions. 
 Typical conditions on the solid boundary are the no-slip condition for  velocity   and  the no-spin condition for microrotation,  which read
$$
\hat u=0\quad\hbox{and}\quad \hat w=0,
$$
which respectively imply that the fluid adheres to the solid boundary and that the fluid elements can not rotate on the fluid-solid interface. However, a  more general boundary condition for  microrotation was introduced to take into account the rotation of the microelements on the solid boundary, which is
effectively proved to be in good accordance with experiments, see Bessonov \cite{Bessonov, Bessonov2} and Migun {\it et al.} \cite{Migun, Migun2}.  This condition, called non-zero spin condition, reads
 \begin{equation}\label{non-zero-spin}
 \hat w\times n={\alpha\over 2}\rot(\hat u)\times n,\quad \hat w\cdot n=0,
 \end{equation}
 where $n$ is a normal unit vector to the boundary and the coefficient $\alpha$ describes the interaction between the given fluid and solid. 
 
It should be noted that in  the previous studies the
no-slip condition for  velocity on the solid surface combined with non-zero spin condition (\ref{non-zero-spin}) for  microrotation is assumed. However slippage is experimentally observed in various systems at fluid-solid interfaces and can strongly influence hydrodynamic behavior in microscale and nanoscale flows. For this reason, no-slip condition for  velocity should be replaced by a more general relation.  In this sense, several boundary conditions have been considered to model the observed slippage, most of them include limited yield stress or retain slippage value proportional to the shear stress.  But there is a new interpretation of the  observed slippage in micropolar fluids, expressed in terms
of the microrotation, by introducing  a new slippage condition for velocity compatible with  non-zero spin boundary condition  for microrotation. This condition was proposed   in  Bayada {\it et al.} \cite{Bayada_NewModels, Bayada_NewModels2} in the framework of lubrication,  and  allows a slippage in the tangential
direction and  retains a non-penetration condition in the normal
direction $n$ ($\delta$ is a real parameter)
 \begin{equation}\label{new_bc_intro}
 \hat u\times n=\delta \rot(\hat w)\times n,\quad \hat u\cdot n=0.
 \end{equation}

On the other hand, the behavior of   fluid flows in porous media is of great importance in industrial and engineering applications. As is well known, classical (Newtonian) Darcy's law   is generally considered for modelling of flow through a porous media, see Darcy \cite{Darcy}. By using homogenization techniques, the mathematical derivation of such Darcy's law was obtained in  Tartar \cite{Tartar} assuming no-slip boundary conditions
$$\begin{array}{l}
\hat u=K(\hat f-\nabla \hat p),\quad {\rm div}(\hat u)=0,
\end{array}
$$
where the matrix coefficient $K$ is calculated by using Newtonian local problems. In addition, problems with different types of  slippage conditions for Newtonian fluids in porous media have been studied by several authors giving rise to a wide range of Darcy's laws. More precisely, the Navier-Stokes (or Stokes)  
flow in a periodic porous media with Fourier boundary conditions on the boundary of the
obstacles was studied in  Conca \cite{Conca} by using the method of oscillating test functions and two-scale method with asymptotic expansion of the solution. The case of classical slip boundary conditions   was treated in Allaire \cite{Allaire} by means of the method of oscillating test functions.   The case of non-homogeneous slip boundary conditions was considered in Cioranescu {\it et al.} \cite{Ciora_Don_Ene}  combining the method of oscillating test functions with the technique introduced in Vanninathan \cite{Vaninathan} to treat the surface integrals. Finally, the case of non-homogeneous slip boundary conditions was revisited  in   Capatina and Ene  \cite{Capatina_Ene} and  Zaki   \cite{Zaki}  by using the periodic unfolding method together with the boundary unfolding operator which allows to treat quite elementary  the surface integrals, see Cioranescu {\it et al.} \cite{Ciora_donato_Griso, Ciora_Don_Zak, Cioranescu_book}.

Although the behavior of micropolar fluid flows in porous media become of great practical relevance,  the literature on the modelling of such type of problem by homogenization methods is far less complete. Lukaszewicz \cite{Luka} rigorously derived
the following version of the classical Darcy law  by using the  two-scale convergence method
$$\begin{array}{l}
\hat u=K^{(1)}(\hat f-\nabla \hat p)+K^{(2)}\hat g,\quad {\rm div}(\hat u)=0,\quad 
 \hat w=L^{(1)}(\hat f-\nabla \hat p)+L^{(2)}\hat g,\\
\end{array}
$$
where  the matrix coefficients $K^{(k)}$ and $L^{(k)}$, $k=1, 2$, are calculated by using micropolar local problems. We also refer to Aganovic and Tutek \cite{Aganovic} for the nonstationary case  and  to Bayada {\it et al.} \cite{Bayada_Cham_Gamo2} for the micropolar effects in the coupling of a thin film past a porous media.

Previous studies obtained  different Darcy's laws for micropolar fluids by assuming on the obstacles of the porous media the no-slip conditon for velocity and the no-spin condition for microrotation, not allowing to capture  the microscopic behavior of the fluid near the boundary of the obstacles. Thus, the goal of this paper is first to establish existence and uniqueness of solution of the micropolar system in the considered porous media  by assuming  non-zero spin boundary condition (\ref{non-zero-spin}) and  new slippage condition (\ref{new_bc_intro}) on the boundary of the obstacles, and then to derive  a generalized Darcy's law by means of a combination of the periodic unfolding method with the boundary unfolding operator to treat the surface integrals.  As far as the author knows, this is the first attempt to carry out such an homogenization analysis for micropolar fluids in porous media, which is the main novelty of the work, and could be instrumental for understanding the effects on this type of non-Newtonian fluid flows taking into account the boundary of the obstacles. 

 The structure of the paper is as follows. In Section 2, we make an introduction of the problem and its setting. In Section 3,  we give 
 the main results of the paper, i.e. the existence and uniqueness of solution (Theorem \ref{lema_coercive}) and the asymptotic behavior of the solution (Theorem \ref{thm_darcy_law}).  The proof of the corresponding results are given in Section 4. The paper ends with a list of references.
 
\section{The setting of the problem}
\paragraph{Definition of the domain.}Let $\Omega$ be a bounded connected open set in $\mathbb{R}^3$, with smooth enough boundary $\partial\Omega$. Denote $Y=(0,1)^3$ and  $F$ an open connected subset of $Y$ with a $C^{1,1}$ boundary
 $\partial F$, such that $\bar F\subset Y$. We denote $Y^*=Y\setminus \bar F$.

For $k\in\mathbb{Z}^2$, each cell $Y_{k,\varepsilon}=\varepsilon k+\varepsilon Y$ is similar to the unit cell $Y$ rescaled to size $\varepsilon$ and $F_{k,\varepsilon}=\varepsilon k+\varepsilon F$ is similar to $F$ rescaled to size $\varepsilon$. We denote $Y^*_{k,\varepsilon}=Y_{k,\varepsilon}\setminus \bar F_{k,\varepsilon}$. 

We denote by $\tau(\bar F_{k,\ep})$ the set of all translated images of $\bar F_{k,\ep}$. The set $\tau(\bar F_{k,\ep})$ represents the obstacles
 in $\mathbb{R}^3$. The porous media is defined by $\Omega_\ep=\Omega\setminus \bigcup_{k\in \mathcal{K}_\ep}\bar F_{k,\varepsilon}$, where   $\mathcal{K}_\ep:=\left\{k\in\mathbb{Z}^3\,:\, Y_{k,\varepsilon}\cap \Omega\neq \emptyset\right\}$. By this construction, $\Omega_\varepsilon$ is a periodically perforated domain with obstacles of the same size as the period.

 We make the assumption that the obstacles $\tau(\bar F_{k,\ep})$ do no intersect the boundary $\partial\Omega$. We denote by $F_\varepsilon$ the set of all the obstacles contained in $\Omega_\ep$. Then,  $F_\varepsilon=\cup_{k\in \mathcal{K}_\ep}\bar F_{k,\ep}$.

 We define $n$ the outside normal vector to $\partial F$. We denote by  $n_\ep(x)=n(x/\ep)$ the outside normal vector (extended by periodicity) to $\partial F_\ep$.

\paragraph{Statement of the problem.} We consider that the micropolar fluid flow is described by the following linearized micropolar equations
in $\Omega_\varepsilon$, taking into account the dependence of $\varepsilon$,
\begin{eqnarray}
-(\nu+\nu_r)\Delta \hat u_\ep+\nabla \hat p_\ep=2\nu_r{\rm rot}(\hat w_\ep)+\hat  f_\ep &\quad\hbox{in}\quad\Omega_\ep,\label{system_1_dimension_1}\\
\noame
{\rm div}(\hat u_\ep)=0&\quad\hbox{in}\quad\Omega_\ep,\label{system_1_dimension_2}\\
\noame
-(c_a+c_d)\Delta\hat  w_\ep+4\nu_r\hat w_\ep=2\nu_r{\rm rot}(\hat u_\ep) + \hat g_\ep&\quad\hbox{in}\quad\Omega_\ep.\label{system_1_dimension_3}
\end{eqnarray}
As discussed in the introduction, we impose the following boundary conditions for velocity and microrotation on the surface of the obstacles
\begin{eqnarray}
\displaystyle \hat w_\varepsilon\times n_\ep={\alpha\over 2} \rot(\hat u_\varepsilon)\times n_\ep & \hbox{ on }\partial F_\varepsilon,\label{BC_holes_dimension_1}\\
\noame
\rot(\hat w_\varepsilon)\times n_\ep=\displaystyle{2\nu_r\over c_a+c_d} \beta (\hat u_\varepsilon\times n_\ep)& \hbox{ on }\partial F_\varepsilon,\label{BC_holes_dimension_2}\\
\noame
\hat u_\varepsilon\cdot n_\ep=0 & \hbox{ on }\partial F_\varepsilon,\label{BC_holes_dimension_3}\\
\noame
\hat w_\varepsilon\cdot n_\ep=0 & \hbox{ on }\partial F_\varepsilon,\label{BC_holes_dimension_4}
\end{eqnarray}
 and the homogeneous boundary conditions on the exterior boundary
\begin{equation}\label{BC_exterior_dimension}
\begin{array}{rl}
\hat u_\varepsilon=0,\quad  \hat w_\varepsilon=0& \hbox{ on }\partial\Omega.
\end{array}
\end{equation}
Notice that the usual no-slip and no-spin boundary conditions for the velocity and microrotation are prescribed in the exterior boundary, while  non-zero spin and new slip boundary conditions are imposed on the boundary of the obstacles. The coefficient $\alpha>0$ appearing in (\ref{BC_holes_dimension_1}) describes the interaction
between the given fluid and solid, it characterizes microrotation retardation on the solid surfaces. In \cite{Bessonov} it was proposed to connect it with the different viscosity coefficients, which allows to give a certain physical sense and to determine the real limits of its value. The coefficient $\beta>0$ in (\ref{BC_holes_dimension_2}) is a characteristic of a slippage and allows the control of the slippage at the boundary of the obstacles when the value $\hat u_\ep$ is not zero.

\paragraph{Mathematical justification of the new slip boundary condition (\ref{BC_holes_dimension_2}).} By assuming condition (\ref{BC_holes_dimension_1}), the supplementary condition (\ref{BC_holes_dimension_2}) on the boundary of the obstacles is needed to close the system.  The derivation of such boundary condition follows arguments from \cite{Bayada_NewModels} which is given in the context of lubrication by applying the non-zero spin condition to a flat surface. The idea is to consider $\psi\in H^1(\Omega_\varepsilon)^3$ , $\psi=0$ on $\partial\Omega$ and $\psi\cdot n_\ep=0$ on $\partial F_\varepsilon$ and recall the following identities
\begin{equation}\label{prop_1}
-\Delta \varphi={\rot}(\rot(\varphi))-\nabla\,{\rm div}(\varphi)\quad\forall\, \varphi\in \mathcal{D}(\Omega_\varepsilon)^3,
\end{equation}
and 
\begin{equation}\label{prop_1_2}{\rm div}(\varphi\times \psi)=\psi\cdot \rot (\varphi)-\varphi\cdot \rot(\psi)\,.
\end{equation}
Integrating by parts and taking into account the divergence theorem, we have 
$$\int_{\Omega_\varepsilon}{\rm div}(\varphi\times \psi)\,dx=\int_{\partial\Omega}(\varphi\times \psi)\cdot n_\ep\,d\sigma-\int_{\partial F_\varepsilon}(\varphi\times \psi)\cdot n_\ep\,d\sigma=\int_{\partial F_\varepsilon}(\varphi\times n_\ep)\cdot \psi\,d\sigma,$$
and taking into account the last identity and (\ref{prop_1_2}), we get
\begin{equation}\label{prop_2}
\int_{\Omega_\varepsilon}\rot(\varphi)\cdot\psi\,dx=\int_{\Omega_\varepsilon}\rot(\psi)\cdot\varphi\, dx +\int_{\partial F_\varepsilon}(\varphi\times n_\ep)\cdot \psi\,d\sigma\quad \forall\, (\varphi,\psi)\in H^1(\Omega_\varepsilon)^3\times H^1(\Omega_\varepsilon)^3.
\end{equation}
Thus, multiplying (\ref{system_1_dimension_3}) by  test function $\psi$  and using identities (\ref{prop_1}) and (\ref{prop_2}), we get
$$
\begin{array}{l}\displaystyle
(c_a+c_d)\left(\int_{\Omega_\varepsilon}\rot(\hat w_\varepsilon)\cdot \rot(\psi)\,dx+\int_{\Omega_\varepsilon}{\rm div}(\hat w_\varepsilon)\cdot {\rm div}(\psi)\,dx
+\int_{\partial F_\varepsilon}\left(\rot(\hat w_\varepsilon)\times n_\ep\right)\cdot \psi\,d\sigma\right)\\
\noame
\displaystyle+4\nu_r\int_{\Omega_\varepsilon}\hat w_\varepsilon\cdot \psi\,dx=2\nu_r\left(\int_{\Omega_\varepsilon} \rot(\psi)\cdot \hat u_\varepsilon\,dx
+\int_{\partial F_\varepsilon}(\hat u_\varepsilon\times n_\ep)\cdot \psi\,d\sigma\right)+\int_{\Omega_\ep}\hat g_\ep\cdot \psi\,dx
\end{array}
$$
In this equation, the unknown terms $\rot(\hat w_\varepsilon)$  on $\partial F_\varepsilon$ prevent a well-posed variational formulation being obtained. It is then possible to cancel the boundary terms on $\partial F_\varepsilon$ by assuming
$$
\rot(\hat w_\varepsilon)\times n_\ep={2\nu_r\over c_a+c_d}(\hat u_\varepsilon\times n_\ep).
$$
Finally, similarly to \cite{Bayada_NewModels}, we assume the slippage condition  (\ref{BC_holes_dimension_2})  on the boundary of the obstacles $\partial F_\varepsilon$ by introducing an additional parameter $\beta>0$ which will enable the influence of this new condition to be controlled when the value $\hat u_\ep$ is not zero.

\paragraph{Dimensionless equations.} It has been observed (see e.g. \cite{Bayada_Cham_Gamo, Bayada_Luka}) that the magnitude
of the viscosity coefficients appearing in equations (\ref{system_1_dimension_1})-(\ref{system_1_dimension_3}) may influence
the effective flow. Thus, it is reasonable to work with the system
written in a non-dimensional form. In view of that, we introduce
$$\begin{array}{c}
\displaystyle u_\ep={\hat u_\ep},\quad p_\ep={\hat p_\ep\over \nu+\nu_r},\quad w_\ep={\hat w_\ep},\quad f_\ep={\hat f_\ep\over \nu+\nu_r},\quad g_\ep={\hat g_\ep\over \nu+\nu_r}\,,\quad
\displaystyle N^2={\nu_r\over \nu+\nu_r},\quad R_M={c_a+c_d\over \nu+\nu_r}.$$
\end{array}
$$
Dimensionless (non-Newtonian) parameter $N^2$ characterizes the coupling between the equations for the velocity and microrotation
and it is of order $\mathcal{O}(1)$, in fact $N^2$ lies between zero and one. The second dimensionless parameter, denoted
by $R_M$ is  related to the characteristic length of the microrotation effects and is compared with small
parameter $\ep$. Thus, we assume that $R_M =\mathcal{O}(\ep^2)$, namely
\begin{equation}\label{RM}
R_M=\ep^2 R_c,\quad\hbox{with }R_c=\mathcal{O}(1).
\end{equation}
This case is the situation that is commonly introduced to study the micropolar fuid because the angular momentum equation  shows a strong coupling between velocity and microrotation in the limit, see \cite{Aganovic, Luka}.
\\

The flow equations (\ref{system_1_dimension_1})-(\ref{system_1_dimension_3}) now have the following form
\begin{eqnarray}
-\Delta u_\ep+\nabla p_\ep=2N^2{\rm rot}(w_\ep)+ f_\ep &\quad\hbox{in}\quad\Omega_\ep,\label{system_1_1}\\
\noame
{\rm div}(u_\ep)=0&\quad\hbox{in}\quad\Omega_\ep,\label{system_1_2}\\
\noame
-\ep^2R_c\Delta w_\ep+4N^2w_\ep=2N^2{\rm rot}(u_\ep) + g_\ep&\quad\hbox{in}\quad\Omega_\ep.\label{system_1_3}
\end{eqnarray}
Concerning the body force and body torque, in order to obtain appropriate estimates, given $f,g\in L^2(\Omega)^3$, we make the following assumptions
\begin{equation}\label{body_forces}
f_\ep(x)=\ep^{-1}f(x),\quad g_\ep(x)= g(x),\quad \hbox{a.e. }x\in \Omega_\ep\,.
\end{equation} 
The corresponding boundary conditions on the boundary of the obstacles read
\begin{eqnarray}
\displaystyle w_\varepsilon\times n_\ep={\alpha\over 2} \rot(u_\varepsilon)\times n_\ep & \hbox{ on }\partial F_\varepsilon,\label{BC_holes_1}\\
\noame
\rot(w_\varepsilon)\times n_\ep=\displaystyle{2N^2\over \ep^2R_c} \beta (u_\varepsilon\times n_\ep)& \hbox{ on }\partial F_\varepsilon,\label{BC_holes_2}\\
\noame
u_\varepsilon\cdot n_\ep=0 & \hbox{ on }\partial F_\varepsilon,\label{BC_holes_3}\\
\noame
w_\varepsilon\cdot n_\ep=0 & \hbox{ on }\partial F_\varepsilon,\label{BC_holes_4}
\end{eqnarray}
and the boundary conditions on the exterior boundary read as follows
\begin{equation}\label{BC_exterior}
\begin{array}{rl}
u_\varepsilon=0,\quad w_\varepsilon=0& \hbox{ on }\partial\Omega.
\end{array}
\end{equation}

\paragraph{Functional setting.}Due to the boundary conditions (\ref{BC_holes_3}) and (\ref{BC_holes_4}), we introduce the functional spaces $V_\ep$ and  $V_\ep^0$ given by
$$\begin{array}{l}
V_\varepsilon=\left\{\varphi\in H^1(\Omega_\varepsilon)^3\,:\, \varphi=0\ \hbox{ on }\partial\Omega,\quad \varphi\cdot n_\ep=0\ \hbox{ on } \partial F_\varepsilon\right\}\,,\quad 
V_\varepsilon^0=\left\{\varphi\in V_\varepsilon\,:\, {\rm div}(\varphi)=0\right\}\,,
\end{array}$$
equipped with the norm $\|\nabla\varphi\|_{L^2(\Omega_\varepsilon)^{3\times 3}}$ and the $L^2_0$ the space of functions of $L^2$ with null integral equipped with the norm of $L^2$.   Let $C^\infty_{\rm per}(Y^*)$ be the space of infinitely differentiable functions in $\mathbb{R}^3$ that are $Y$-periodic. By $L^2_{\rm per}(Y^*)$ (resp. $H^1_{\rm per}(Y^*)$) we denote its completion in the norm $L^2(Y^*)$ (resp. $H^1(Y^*)$). We denote by $L^2_{0,{\rm per}}(Y^*)$  the space of functions in $L^2_{\rm per}(Y^*)$ with null integral. We also define the spaces $V_Y$ and $V_Y^0$ are given by
$$V_Y=\left\{\varphi\in H^1_{\rm per}(Y^*)^3\,:\,  \varphi\cdot n=0\hbox{ on }\partial F\right\},\quad V_Y^0=\left\{\varphi\in V_Y\,:\,  {\rm div}_y(\varphi)=0\hbox{ in }Y^*\right\}.$$

\section{Main results} \label{sec:MainResult}
In this section we give the main results of this paper. First, the existence and uniqueness of solution of  problem  (\ref{system_1_1})-(\ref{BC_exterior})  is established in Theorem \ref{lema_coercive} and then,  the homogenized model is given in Theorem \ref{thm_darcy_law}. The proof of the corresponding
results are given in the next section.

In order to prove, for each value of $\ep$,  the existence and uniqueness of solution of problem (\ref{system_1_1})-(\ref{BC_exterior}), instead
of working directly with the classical variational formulation, we will work with an equivalent variational formulation. 

\begin{proposition}\label{prop_form_var} Sufficiently  regular solutions of (\ref{system_1_1})-(\ref{BC_exterior}) satisfy the following weak formulation: 

Find $(u_\ep, w_\ep, p_\ep)\in V_\ep^0\times V_\ep\times L^2_0(\Omega_\ep)$ such that 
\begin{equation}\label{var_form_v}
\begin{array}{l}\displaystyle
\int_{\Omega_\varepsilon}\rot(u_\ep)\cdot \rot(\varphi)\,dx- \int_{\Omega_\ep}p_\ep{\rm div}(\varphi)\,dx-2N^2\int_{\Omega_\varepsilon}\rot(\varphi)\cdot w_\ep\,dx
\\
\noame\displaystyle
\qquad
+2\left({1\over \alpha}-N^2\right)\int_{\partial F_\varepsilon}(w_\ep\times n_\ep)\cdot\varphi\,d\sigma= \ep^{-1}\int_{\Omega_\ep}f\cdot \varphi\,dx\,,\quad \forall \varphi\in V_\ep\,,
\end{array}
\end{equation}
\begin{equation}\label{var_form_w}
\begin{array}{l}\displaystyle
\ep^2R_c\int_{\Omega_\varepsilon}\rot(w_\ep)\cdot \rot(\psi)\,dx+\ep^2R_c\int_{\Omega_\varepsilon}{\rm div}(w_\ep)\cdot {\rm div}(\psi)\,dx+4N^2\int_{\Omega_\varepsilon}w_\ep\cdot\psi\,dx\\
\noame\displaystyle
\qquad-2N^2\int_{\Omega_\varepsilon}\rot(\psi)\cdot u_\ep\,dx + 2N^2(\beta-1)\int_{\partial F_\varepsilon}(u_\ep\times n_\ep)\cdot \psi\,d\sigma
= \int_{\Omega_\ep}g\cdot \psi\,dx\,,\quad \forall \psi\in V_\ep\,.
\end{array}
\end{equation}
\end{proposition}

We give the result of the existence and uniqueness of solution of problem (\ref{var_form_v})-(\ref{var_form_w}). 

\begin{theorem}\label{lema_coercive}
 Let $\gamma={1\over \alpha}-N^2-N^2\beta$. Assume that the asymptotic regimes (\ref{RM}) and (\ref{body_forces}) hold. Then, for any $\alpha$ and $\beta$  satisfying
 \begin{equation}\label{condition_existence}
\gamma^2<{R_c(1-N^2)\over K^2}\,,
 \end{equation}
there exists a unique solution $(u_\varepsilon, w_\varepsilon, p_\varepsilon)\in V_\ep^0\times V_\ep\times L^2_0(\Omega_\ep)$ of problem (\ref{var_form_v})-(\ref{var_form_w}), with $K=C_{pt}C_g$ where $C_{pt}$ and $C_g$ are the trace and the Gaffney constants  given, respectively,   in Corollary \ref{cortrace} and Lemma \ref{Lemma:Gaffney}, both placed in the next section.
 \end{theorem}

We describe the asymptotic behavior of the sequences $ u_\ep$, $w_\ep$ and $ p_\ep$ when $\ep$ tends to zero.  To do this, we take into account that the sequence of solutions $(u_\ep, w_\ep,  p_\ep)\in V_\ep^0\times V_\ep\times L^2_0(\Omega_\ep)$ is not defined in a fixed domain independent of $\ep$ but rather in a varying set $\Omega_\ep$. Thus, in order to pass to the limit when $\ep$ tends to zero, convergences in fixed Sovolev spaces (defined in $\Omega$) are used, which require first that $(u_\ep, w_\ep,  p_\ep)$ be extended to the whole domain $\Omega$. Then, we define an extension $(U_\ep, W_\ep,  P_\ep)$ of $(u_\ep, w_\ep, p_\ep)$ on $\Omega$ which coincides with $(u_\ep, w_\ep,   p_\ep)$ on $\Omega_\ep$.

\begin{theorem}\label{thm_darcy_law} Assume that the asymptotic regimes (\ref{RM}) and (\ref{body_forces}) and condition (\ref{condition_existence}) hold. Then, the whole sequences of extensions $(\ep^{-1}U_\ep, W_\ep)$  and $\ep P_\ep$ of the solution of problem (\ref{system_1_1})-(\ref{BC_exterior}) converge weakly to $(u, w)$ in $L^2(\Omega)^3\times L^2(\Omega)^3$ and strongly to $p$ in $L^2(\Omega)$ respectively. Moreover, it holds
\begin{equation}\label{Darcy_law_vel_mic}
u(x)=K^{(1)}\left(f(x)-\nabla p(x)\right)+ K^{(2)} g(x),\quad w(x)=L^{(1)}\left(f(x)-\nabla p(x)\right)+ L^{(2)} g(x)\quad\hbox{ in }\Omega,
\end{equation}
and also, $p\in H^1(\Omega)\cap L^2_0(\Omega)$ is the unique solution of the Darcy equation
\begin{equation}\label{Darcy_law}
\left\{\begin{array}{rl}
{\rm div}\left(K^{(1)}\left(f(x)-\nabla p(x)\right)+ K^{(2)} g(x)\right)=0 & \hbox{in }\Omega,\\
\noame
 \left( K^{(1)}\left(f(x)-\nabla p(x)\right)+ K^{(2)} g(x)\right)\cdot n=0& \hbox{on }\partial \Omega.
\end{array}\right.
\end{equation}
The matrix coefficients $K^{(k)}$, $L^{(k)}\in \mathbb{R}^{3\times 3}$, $k=1,2$, where $K^{(1)}$ is positive definite, are given by
$$K^{(k)}_{ij}=\int_{Y^*}u^{i,k}_j(y)\,dy,\quad L^{(k)}_{ij}=\int_{Y^*}w^{i,k}_j(y)\,dy,\quad i,j=1,2,3,$$
with $(u^{i,k},w^{i,k}, \pi^{i,k})\in V^0_Y\times V_Y\times L^2_{0,{\rm per}}(Y^*)$, $k=1,2$,  $i=1,2,3$, the unique solution of local micropolar problem 
\begin{equation}\label{local_problems}
\left\{\begin{array}{rl}
-\Delta_y u^{i,k}+\nabla_y \pi^{i,k}-2N^2 \rot_y(w^{i,k})=e_i\delta_{1k}&\hbox{in }Y^*,\\
\noame
-R_c\Delta_y w^{i,k}+4N^2 w^{i,k}-2N^2 \rot_y(u^{i,k})=e_i\delta_{2k}&\hbox{in }Y^*,\\
\noame
{\rm div}_y(u^{i,k})=0&\hbox{in }Y^*,\\
\noame
\displaystyle w^{i,k}\times n={\alpha\over 2} \rot(u^{i,k})\times  n& \hbox{ on }\partial F,\\
\noame
\rot(w^{i,k})\times n=\displaystyle{2N^2\over R_c} \beta (u^{i,k}\times n)& \hbox{ on }\partial F,\\
u^{i,k}\cdot n=0 & \hbox{ on }\partial F,\\
\noame
w^{i,k}\cdot n=0 & \hbox{ on }\partial F.
\end{array}\right.
\end{equation}
\end{theorem}

\section{Proof of the main results}\label{sec:proofs}
\begin{proof}[Proof of Proposition \ref{prop_form_var}]
First, from (\ref{system_1_2}) and the boundary conditions (\ref{BC_holes_3}) and (\ref{BC_holes_4}),  solutions of (\ref{system_1_1})-(\ref{BC_exterior}) are in $V_\ep^0\times V_\ep$. 

Next, to obtain (\ref{var_form_v}), we take $\varphi\in V_\ep$ as test function in (\ref{system_1_1}) and using (\ref{prop_1}) and (\ref{prop_2}) we get 
$$\begin{array}{l}\displaystyle
\int_{\Omega_\varepsilon}\rot(u_\varepsilon)\cdot \rot(\varphi)\,dx- \int_{\Omega_\ep}p_\ep{\rm div}(\varphi)\,dx +\int_{\partial F_\varepsilon}\left(\rot(u_\varepsilon)\times n_\ep\right)\cdot \varphi\,d\sigma\\
\noame
\displaystyle-2N^2\int_{\Omega_\varepsilon} \rot(\varphi)\cdot w_\varepsilon\,dx
-2N^2\int_{\partial F_\varepsilon}(w_\varepsilon\times n_\ep)\cdot \varphi\,d\sigma= \ep^{-1}\int_{\Omega_\ep}f\cdot \varphi\,dx\,.
\end{array}$$
Thus, taking into account the boundary condition (\ref{BC_holes_1}), we derive equation (\ref{var_form_v}).

Finally, to obtain (\ref{var_form_w}), we take $\psi\in V_\ep$ as test function in (\ref{system_1_3}) and proceeding as above we get 
$$\begin{array}{l}\displaystyle
\ep^2R_c\int_{\Omega_\varepsilon}\rot(w_\varepsilon)\cdot \rot(\psi)\,dx+\ep^2R_c\int_{\Omega_\varepsilon}{\rm div}(w_\varepsilon)\, {\rm div}(\psi)\,dx
+\ep^2R_c\int_{\partial F_\varepsilon}\left(\rot(w_\varepsilon)\times n_\ep\right)\cdot \psi\,d\sigma\\
\noame
\displaystyle+4N^2\int_{\Omega_\varepsilon}w_\varepsilon\cdot \psi\,dx-2N^2\int_{\Omega_\varepsilon} \rot(\psi)\cdot u_\varepsilon\,dx
-2N^2\int_{\partial F_\varepsilon}(u_\varepsilon\times n_\ep)\cdot \psi\,d\sigma=\int_{\Omega_\ep}g\cdot \psi\,dx\,.
\end{array}$$
Taking into account the boundary condition (\ref{BC_holes_2}), we derive equation (\ref{var_form_w}).
\end{proof}

  Before proving the result concerning existence and uniqueness of solution,  we give several technical lemmas which will also be used to obtain {\it a priori} estimates of the solution. First, we recall a result about a trace result on the boundary of the obstacles $\partial F_\ep$ whose proof can be found in \cite{Allaire, Anguiano, Anguiano2, Cioranescu_book}.

\begin{lemma}
There exists a positive constant $C_t$ independent of $\ep$, such that for every $v\in H^1(\Omega_\ep)^3$,
\begin{equation}\label{estimate_buo}
\|v\|^2_{L^2(\partial F_\ep)^3}\leq \ep^{-1}C_t\left(\|v\|^2_{L^2(\Omega_\ep)^3}+\ep^2\|D v\|^2_{L^2(\Omega_\ep)^{3\times 3}}\right).
\end{equation}
\end{lemma}

We note that the fact that the normal component of the function is equal zero on $\partial F_\ep$ is not used in the previous estimate, so it holds true in a more general context. However, the next inequalities in $\Omega_\ep$  makes use of this condition and will be also  used for obtaining a priori estimates for the velocity and microrotation. Thus,  we first recall the version of  Poincar\'e's inequality given in \cite{Allaire, Capatina_Ene}. 

\begin{lemma} There exists a positive constant $C_p$, independent of $\ep$, such that for every $v \in V_\ep$,
\begin{equation}\label{Poincare}
\|v\|_{L^2(\Omega_\ep)^3}\leq \ep C_p\|D v\|_{L^2(\Omega_\ep)^{3\times 3}}.
\end{equation}
\end{lemma}

As a consequence of  previous results, we deduce the following result.
\begin{corollary} \label{cortrace}For every $v\in V_\ep$, the following estimate holds 
\begin{equation}\label{Unfolding_boundary_1}
\|v\|_{L^2(\partial F_\ep)^3}\leq (\ep C_{pt})^{1\over 2}\|D v\|_{L^2(\Omega_\ep)^{3\times 3}}, 
\end{equation}
where the positive constant $C_{pt}=C_t^2(C_p^2+1)^2$ with $C_t$ given in (\ref{estimate_buo}) and $C_p$ given in (\ref{Poincare}).
\end{corollary}

We also give an estimate of the derivative in terms of the divergence and the rotational, necessary to prove the coercivity of the variational formulation. It has different names in the literature, e.g.  Gaffney's, Maxwell's or Friedrichs' inequality, see \cite{Amorouche, Boyer_Fabrie, Duvaut_Lions, Grisvard}. 
\begin{lemma} \label{Lemma:Gaffney}There exists a positive constant $C_g$, independent of $\ep$, such that  for every $v \in V_\ep$,
\begin{equation}\label{Maxwell_inequality}
\|D v\|^2_{L^2(\Omega_\ep)^{3\times 3}}\leq C_g\left( \|{\rm div}(v)\|^2_{L^2(\Omega_\ep)}+ \|\rot(v)\|^2_{L^2(\Omega_\ep)^3}\right).
\end{equation}
Moreover,  for every $v \in V_\ep^0$ it holds
\begin{equation}\label{Maxwell_inequality_div0_cor}
\|D v\|_{L^2(\Omega_\ep)^{3\times 3}}\leq C_g \|\rot(v)\|_{L^2(\Omega_\ep)^3}.
\end{equation}
\end{lemma}

\begin{proof} For every function $v\in H^1(Y^*_{k,1})^3$ such that $v\cdot n=0$ on $\partial F_{k,1}$, using Theorem IV.4.8. in \cite{Boyer_Fabrie} (see also Chapter 7, Lema 6.1 in \cite{Duvaut_Lions}), the fact that $Y^*_{k,1}$ is simply connected and  that the boundary $\partial F_{k,1}$ is $C^{1,1}$, we have, for every $k\in \mathbb{Z}^3$, that
\begin{equation}\label{proof_gaffney}
\|D v\|^2_{L^2(Y^*_{k,1})^{3\times 3}}\leq C_g\left( \|{\rm div}(v)\|^2_{L^2(Y^*_{k,1})}+ \|\rot(v)\|^2_{L^2(Y^*_{k,1})^3}\right)\,,
\end{equation}
where the positive constant $C_g$ only depends on $Y^*_{k,1}$.

By the change of variable
$$
y={x\over \ep},\quad dy=\ep^{-3}dx,\quad \partial_y=\ep\partial_x\,,
$$
we rescale (\ref{proof_gaffney}) from $Y^*_{k,1}$ to $Y^*_{k,\ep}$ and from $F_{k,1}$ to $F_{k,\ep}$.  This yields
$$
\|D v\|^2_{L^2(Y^*_{k,\ep})^{3\times 3}}\leq C_g\left( \|{\rm div}(v)\|^2_{L^2(Y^*_{k,\ep})}+ \|\rot(v)\|^2_{L^2(Y^*_{k,\ep})^3}\right)\,.
$$
Summing the inequalities, for every $k\in \mathcal{K}_\ep$, gives (\ref{Maxwell_inequality}) (estimate (\ref{Maxwell_inequality_div0_cor}) is straightforward).

In fact, we must consider separately the periods containing a portion of $\partial\Omega$, but they yield at a distance $O(\ep)$ of
$\partial\Omega$, where $v$ is zero, and then the corresponding inequality is immediately obtained.
\end{proof}

\begin{proof}[Proof of Theorem \ref{lema_coercive}] To prove, for each value of $\ep>0$, the existence and uniqueness of a weak solution 
$(u_\ep, w_\ep, p_\ep)$ of problem (\ref{var_form_v})-(\ref{var_form_w}), we shall apply classical results given in \cite{Girault_Raviart}.  To this purpose, we introduce the  following equivalent mixed variational form:

Find $(u_\ep,w_\ep,p_\ep)\in V_\ep\times V_\ep\in L^2_0(\Omega_\ep)$ such that
\begin{equation}\label{bilinear_var_form0}
\begin{array}{rcll}
\mathcal{A}_\ep(u_\ep, w_\ep;\varphi,\psi)+\mathcal{B}_\ep((\varphi,\psi),p_\ep )&=&\mathcal{L}_\ep(\varphi,\psi)& \forall\, (\varphi,\psi)\in V_\ep\times V_\ep,
\\
\noame
\mathcal{B}_\ep((u_\ep,w_\ep),q_\ep)&=&0&    \forall\, q_\ep\in L^2_0(\Omega_\ep)\,,
\end{array}
\end{equation}
where 
\begin{equation}\label{Def_A}\begin{array}{rl}\displaystyle
\mathcal{A}_\ep(u_\ep,w_\ep;\varphi,\psi)=&\displaystyle 
\int_{\Omega_\varepsilon}\rot(u_\ep)\cdot \rot(\varphi)\,dx-2N^2\int_{\Omega_\varepsilon}\rot(\varphi)\cdot w_\ep\,dx\\
\noame&
\displaystyle +\ep^2R_c\int_{\Omega_\varepsilon}\rot(w_\ep)\cdot \rot(\psi)\,dx+\ep^2R_c\int_{\Omega_\varepsilon}{\rm div}(w_\ep)\cdot {\rm div}(\psi)\,dx
\\
\noame& \displaystyle
-2N^2\int_{\Omega_\varepsilon}\rot(\psi)\cdot u_\ep\,dx+4N^2\int_{\Omega_\varepsilon}w_\ep\cdot\psi\,dx\\
\noame &\displaystyle
+2\left({1\over \alpha}-N^2\right)\int_{\partial F_\varepsilon}(w_\ep\times n_\ep)\cdot\varphi\,d\sigma+ 2N^2(\beta-1)\int_{\partial F_\varepsilon}(u_\ep\times n_\ep)\cdot \psi\,d\sigma\,,
\end{array}
\end{equation} and  
$$
\mathcal{B}_\ep((\varphi,\psi),p_\ep)= \displaystyle -\int_{\Omega_\ep}p_\ep\,{\rm div}(\varphi)\,dx,\quad 
\mathcal{L}_\ep(\varphi,\psi)=\displaystyle \ep^{-1}\int_{\Omega_\ep}f\cdot \varphi\,dx+
\int_{\Omega_\ep}g\cdot \psi\,dx\,.
$$
It is easy to prove that the bilinear form $\mathcal{A}_\ep$, $\mathcal{B}_\ep$ and $\mathcal{L}_\ep$  are continuous bilinear forms on $(V_\ep\times V_\ep)^2$, $V_\ep\times V_\ep\times L^2_0(\Omega_\ep)$  and $V_\ep\times V_\ep$ respectively. Denoting
$$\|(\varphi,\psi)\|_{V_\ep^0\times V_\ep}=(\|D \varphi\|^2_{L^2(\Omega_\ep)^{3\times 3}}+\|D \psi\|^2_{L^2(\Omega_\ep)^{3\times 3}})^{1\over 2}.$$ 
classical existence and uniqueness conditions for such a problem given in Theorem 4.1  \cite{Girault_Raviart} are the coerciveness of the form $\mathcal{A}_\ep$ on the subspace $(V_\ep^0\times V_\ep)^2$ and the inf-sup condition.

First, we prove that  $\mathcal{A}_\ep$ is coercive   on $(V_\ep^0\times V_\ep)^2$, i.e. that there exists  $\eta>0$ such that 
\begin{equation}\label{A_coercive_enunciado}\mathcal{A}_\ep(\varphi,\psi;\varphi,\psi)\geq \eta\|(\varphi,\psi)\|^2_{V_\ep^0\times V_\ep}.
\end{equation}

To do this, let us  derive another equivalent expression for $\mathcal{A}_\ep$. By using (\ref{prop_2}) applied to $\int_{\Omega_\ep}\rot(\psi)\cdot u_\ep\,dx$ and using that $\int_{\partial F_\ep}(\psi\times n_\ep)\cdot u_\ep\,d\sigma=-\int_{\partial F_\ep}(u_\ep\times n_\ep)\cdot \psi\,d\sigma$, we have that  $\mathcal{A}_\ep$ defined in (\ref{Def_A}) has the following expression
 \begin{equation}\label{Def_A_2}\begin{array}{rl}\displaystyle
\mathcal{A}_\ep(u_\ep,w_\ep;\varphi,\psi)=&\displaystyle 
\int_{\Omega_\varepsilon}\rot(u_\ep)\cdot \rot(\varphi)\,dx-2N^2\int_{\Omega_\varepsilon}\rot(\varphi)\cdot w_\ep\,dx\\
\noame&
\displaystyle +\ep^2R_c\int_{\Omega_\varepsilon}\rot(w_\ep)\cdot \rot(\psi)\,dx+\ep^2R_c\int_{\Omega_\varepsilon}{\rm div}(w_\ep)\cdot {\rm div}(\psi)\,dx
\\
\noame& \displaystyle
-2N^2\int_{\Omega_\varepsilon}\rot(u_\ep)\cdot \psi\,dx+4N^2\int_{\Omega_\varepsilon}w_\ep\cdot\psi\,dx+2\gamma\int_{\partial F_\varepsilon}(w_\ep\times n_\ep)\cdot\varphi\,d\sigma\,.
\end{array}
\end{equation}
Then, we have
$$
\begin{array}{rl}\displaystyle
\mathcal{A}_\ep(\varphi,\psi;\varphi,\psi)=&\displaystyle 
\int_{\Omega_\varepsilon}|\rot(\varphi)|^2\,dx-4N^2\int_{\Omega_\varepsilon}\rot(\varphi)\cdot \psi\,dx+\ep^2R_c\int_{\Omega_\varepsilon}|\rot(\psi)|^2\,dx\\
\noame&
\displaystyle +\ep^2R_c\int_{\Omega_\varepsilon}|{\rm div}(\psi)|^2\,dx+4N^2\int_{\Omega_\varepsilon}|\psi|^2\,dx+2\gamma\int_{\partial F_\varepsilon}(\psi\times n_\ep)\cdot\varphi\,d\sigma\,.
\end{array}$$
From the Cauchy-Schwartz inequality and (\ref{Unfolding_boundary_1}), we get
$$\int_{\partial F_\varepsilon}(\psi\times n_\ep)\cdot\varphi\,d\sigma\leq \ep C_{pt} \|D\psi\|_{L^2(\Omega_\ep)^{3\times 3}}\|D\varphi\|_{L^2(\Omega_\ep)^{3\times 3}},$$
and by using (\ref{Maxwell_inequality}) and (\ref{Maxwell_inequality_div0_cor}), we deduce
\begin{equation}\label{A_coercive}
\begin{array}{rl}\displaystyle
\mathcal{A}_\ep(\varphi,\psi;\varphi,\psi) \geq & \displaystyle {1\over C_g}\|D \varphi\|_{L^2(\Omega_\ep)^{3\times 3}}^2-4\sqrt{3}N^2\|D\varphi\|_{L^2(\Omega_\ep)^{3\times 3}}\|\psi\|_{L^2(\Omega_\ep)^3}\\
\noame
&\displaystyle +\ep^2{R_c\over C_g}\|D \psi\|_{L^2(\Omega_\ep)^{3\times 3}}^2+4N^2\|\psi\|^2_{L^2(\Omega_\ep)^3}\\
\noame&\displaystyle - 2|\gamma|\ep C_{pt} \|D \varphi\|_{L^2(\Omega_\ep)^{3\times 3}}\|D \psi\|_{L^2(\Omega_\ep)^{3\times 3}}\,.
\end{array}\end{equation}
Now, by condition (\ref{condition_existence}), there exists $c_1>0$ satisfying
$$
{|\gamma|C_{pt} C_g^2\over R_c}<c_1<{1-N^2\over |\gamma|C_{pt}}\,,
$$
and by Young's inequality,
$$\|D \varphi\|_{L^2(\Omega_\ep)^{3\times 3}}\|D \psi\|_{L^2(\Omega_\ep)^{3\times 3}}\leq {c_1\over 2\ep C_g}\|D \varphi\|^2_{L^2(\Omega_\ep)^{3\times 3}}+{\ep C_g\over 2c_1}\|D \psi\|^2_{L^2(\Omega_\ep)^{3\times 3}}\,.$$
Introducing a real number $c_2$ satisfying $0<c_2<\min\{1,1/(3C_g)\}$, and such that
$$
c_1<{1-{N^2\over c_2}\over |\gamma|C_{pt}}\quad \left(\hbox{and so  }{|\gamma|C_{pt} C_g^2\over R_c}<c_1<{1-{N^2\over c_2}\over |\gamma|C_{pt}}<{1-N^2\over |\gamma|C_{pt}}\right),
$$
we also have by Young's inequality
$$\|D \varphi\|_{L^2(\Omega_\ep)^{3\times 3}}\|\psi\|_{L^2(\Omega_\ep)^{3}}\leq {1\over 4\sqrt{3}c_2 C_g}\|D \varphi\|^2_{L^2(\Omega_\ep)^{3\times 3}}+ {\sqrt{3}\,c_2  C_g}\|\psi\|^2_{L^2(\Omega_\ep)^{3}}\,.$$
Going back to estimate (\ref{A_coercive}), we obtain
\begin{equation*}
\mathcal{A}_\ep(\varphi,\psi;\varphi,\psi)\geq A \|D \varphi\|^2_{L^2(\Omega_\ep)^{3\times 3}}+\ep^2 B \|D \psi\|^2_{L^2(\Omega_\ep)^{3\times 3}}+ 4N^2(1-3c_2 C_g)\|\psi\|^2_{L^2(\Omega_\ep)^{3}}\,,
\end{equation*}
where $A$, $B$ are defined as follows
\begin{equation}\label{AB}
A={1\over C_g}\left(1-{N^2\over c_2}-|\gamma|{C_{pt}\, c_1}\right),\quad B={1\over C_g}\left(R_c-{|\gamma|C_{pt}\over c_1}C_g^2\right)\,.
\end{equation}
Thus, we have
\begin{equation}\label{A_coercive_1}
\mathcal{A}_\ep(\varphi,\psi;\varphi,\psi)\geq A \|D \varphi\|^2_{L^2(\Omega_\ep)^{3\times 3}}+\ep^2 B \|D \psi\|^2_{L^2(\Omega_\ep)^{3\times 3}}\,,
\end{equation}
which proves that (\ref{A_coercive_enunciado}) holds.

Next, we prove the inf-sup condition,  i.e. that there exists $\delta>0$ such that 
\begin{equation}\label{infsup}
\inf_{q_\ep\in L^2_0(\Omega_\ep)}\sup_{(u_\ep,\psi)\in V_\ep^2}{\mathcal{B}_\ep((u_\ep,\psi),q_\ep)\over \|(u_\ep,\psi)\|_{V_\varepsilon^2}\|q_\ep\|_{L^2_0(\Omega_\ep)}}\geq \delta.
\end{equation}
Let $(\varphi,\psi)$ belong to $V_\ep\times V_\ep$, and $q_\ep$ to $L^2_0(\Omega_\ep)$, we obtain
$$H_0^1(\Omega_\ep)^3\times\{0\}\subseteq V_\ep\times V_\ep,$$
so that
$$\sup_{(\varphi,\psi)\in V_\ep\times V_\ep}{\int_{\Omega_\ep}{\rm div}\,\varphi\,q_\ep\,dx\over (\|\varphi\|^2_{V_\ep}+\|\varphi\|^2_{V_\ep})^{1\over 2}}
\geq \sup_{(\varphi,0)\in H^1_0(\Omega_\ep)^3\times \{0\}}{\int_{\Omega_\ep}{\rm div}\,\varphi\,q_\ep\,dx\over \|D\varphi\|_{L^2(\Omega_\ep)^{3\times 3}}}
=\sup_{\varphi\in H^1_0(\Omega_\ep)^3}{\int_{\Omega_\ep}{\rm div}\,\varphi\,q_\ep\,dx\over \|D\varphi\|_{L^2(\Omega_\ep)^{3\times 3}}}.$$
According to the inverse of the divergence operator in perforated domains, see for example \cite{Galdi}, for a given $q_\ep\in L^2_0(\Omega_\ep)$, there exists $v_\ep[q_\ep]\in H^1_0(\Omega_\ep)^3$ such that ${\rm div}\,v_\ep[q_\ep]=q_\ep$ in $\Omega_\ep$ and $\|Dv_\ep[q_\ep]\|_{L^2(\Omega_\ep)^{3\times 3}}\leq C\|q_\ep\|_{L^2(\Omega_\ep)}$ for some constant $C>0$ independent of $\ep$. Choosing $\varphi=v_\ep[q_\ep]$, we get 
$$\sup_{(\varphi,\psi)\in V_\ep\times V_\ep}{\int_{\Omega_\ep}{\rm div}\,\varphi\,q_\ep\,dx\over (\|\varphi\|^2_{V_\ep}+\|\varphi\|^2_{V_\ep})^{1\over 2}}
\geq
{\int_{\Omega_\ep}{\rm div}\,v_\ep[q_\ep]\,q_\ep\,dx\over \|Dv_\ep[q_\ep]\|_{L^2(\Omega_\ep)^{3\times 3}}}
={\|q_\ep\|_{L^2(\Omega_\ep)}^2\over \|Dv_\ep[q_\ep]\|_{L^2(\Omega_\ep)^{3\times 3}}}\geq {1\over C}\|q_\ep\|_{L^2(\Omega_\ep)}.$$
This ends the proof.
 \end{proof}

\paragraph{A priori estimates.}\label{sec:estimates}
We establish sharp  {\it a priori} estimates of the solution  in $\Omega_\ep$ and also for extended solution to $\Omega$, which is independent of $\ep$, introducing suitable extension operators.

We give the a priori estimates for velocity and microrotation in $\Omega_\ep$.
\begin{lemma} \label{Estimates_extended_lemma} Assume that the asymptotic regimes (\ref{RM}) and (\ref{body_forces}) and condition (\ref{condition_existence}) hold. Then there exists a positive constant $C$, independent of $\ep$, such that the following estimates for the velocity and microrotation hold
\begin{equation}\label{estimates_u}
\ep^{-1}\|u_\ep\|_{L^2(\Omega_\ep)^{3}}+\|D u_\ep\|_{L^2(\Omega_\ep)^{3\times 3}}\leq C\,,
\end{equation}
\begin{equation}\label{estimates_w}
\|w_\ep\|_{L^2(\Omega_\ep)^{3}}+\ep\|D w_\ep\|_{L^2(\Omega_\ep)^{3\times 3}}\leq C\,.
\end{equation}
\end{lemma}
\begin{proof} To obtain estimates of velocity and microrotation we consider $(\varphi,\psi)=(u_\ep, w_\ep)$ as test functions in the weak formulation
(\ref{var_form_v})-(\ref{var_form_w}) and use (\ref{A_coercive_1}) and Young's inequality to obtain
$$\begin{array}{l}
\displaystyle
A\|D u_\ep\|^2_{L^2(\Omega_\ep)^{3\times 3}}+\ep^2B\|D w_\ep\|^2_{L^2(\Omega_\ep)^{3\times 3}}\\
\noame
\displaystyle\qquad \leq C_p \|f\|_{L^2(\Omega)^3}\|D u_\ep\|_{L^2(\Omega_\ep)^{3\times 3}}+\ep C_p\|g\|_{L^2(\Omega)^3}\|D w_\ep\|_{L^2(\Omega_\ep)^{3\times 3}}\\
\noame
\displaystyle\qquad\leq {C_p^2\over 2A}\|f\|^2_{L^2(\Omega)^3}+{A\over 2}\|D u_\ep\|^2_{L^2(\Omega_\ep)^{3\times 3}}
+{C_p^2\over 2B}\|g\|^2_{L^2(\Omega)^3}+{\ep^2 B\over 2}\|D w_\ep\|^2_{L^2(\Omega_\ep)^{3\times 3}}\,.
\end{array}$$
Then we have,
$$\begin{array}{l}
\displaystyle
{A\over 2}\|D u_\ep\|^2_{L^2(\Omega_\ep)^{3\times 3}}+\ep^2{B\over 2}\|D w_\ep\|^2_{L^2(\Omega_\ep)^{3\times 3}}\leq {C_p^2\over 2A}\|f\|^2_{L^2(\Omega)^3}+{C_p^2\over 2B}\|g\|^2_{L^2(\Omega)^3}\,.
\end{array}$$
Since $A$ and $B$ are bounded,  we deduce than the right hand side of the
previous inequality is bounded by a certain positive constant $C$ independent of $\ep$, which  implies the following estimates in $\Omega_\ep$
$$
\|D u_\ep\|_{L^2(\Omega_\ep)^{3\times 3}}\leq C\,,\quad \|D w_\ep\|_{L^2(\Omega_\ep)^{3\times 3}}\leq C \ep^{-1}\,.
$$
This and Poincar\'e's inequality (\ref{Poincare}) give the estimate for $u_\ep$ and $w_\ep$ respectively. 
\end{proof}

Since the solution $(u_\ep, w_\ep)$ of problem (\ref{system_1_1})-(\ref{BC_exterior}) is defined only in $\Omega_\ep$, we need to extend them to the whole domain $\Omega$. If we had considered the micropolar equations with Dirichlet boundary condition on the obstacles, the velocity and microrotation would be extended by zero in the obstacles. However,  we need another kind of extension
for the case in which the velocity and microrotation are non-zero on the obstacles.  
Thus, we introduce an extension operator which  is classical in the homogenization literature, see \cite{Allaire, Cioranescu_Paulin, Conca,Tartar}.

\begin{lemma}\label{extension_vel}
There exists an extension operator $\Pi_\ep\in \mathcal{L}(H^1(\Omega_\ep)^3;H^1_0(\Omega)^3)$ and a positive constant $C$, independent of $\ep$, such that
$$\begin{array}{c}
\Pi_\ep v(x)=v(x),\quad \hbox{ if }x\in \Omega_\ep,\\
\noame
\|D\Pi_\ep v\|_{L^2(\Omega)^{3\times 3}}\leq C\|D v\|_{L^2(\Omega_\ep)^{3\times 3}},\quad \forall v\in H^1(\Omega_\ep)^3.
\end{array}
$$
\end{lemma}

Taking into account the extension $\Pi_\ep$, we denote by $U_\ep$ the extension $\Pi_\ep u_\ep$ of the velocity $u_\ep$,  and  by $W_\ep$ the extension  $\Pi_\ep w_\ep$ of the microrotation $w_\ep$. Next, we get the following uniform estimates in $\Omega$ as consequence of  Lemmas \ref{Estimates_extended_lemma} and \ref{extension_vel}. 

\begin{corollary}\label{lemma_ext_vel}Assume that the asymptotic regimes (\ref{RM}) and (\ref{body_forces}) and condition (\ref{condition_existence}) hold. Then there exists a positive constant $C$, independent of $\ep$, such that the following estimates for the extensions of velocity and microrotation hold
\begin{equation}\label{estimates_u_ext}
\ep^{-1}\|U_\ep\|_{L^2(\Omega)^{3}}+\|D U_\ep\|_{L^2(\Omega)^{3\times 3}}\leq C\,,
\end{equation}
\begin{equation}\label{estimates_w_ext}
\| W_\ep\|_{L^2(\Omega)^{3}}+\ep \|D W_\ep\|_{L^2(\Omega)^{3\times 3}}\leq {C}\,.
\end{equation}
\end{corollary}

Now, we recall two important results from  \cite{Tartar} which are concerned with the extension of the pressure $p_\ep$ to the whole domain $\Omega$. First, we  define a restriction operator $R_\ep$ from $H^1_0(\Omega)^3$ into $H^1_0(\Omega_\ep)^3$ and then, we extend the gradient of the pressure by duality in $H^{-1}(\Omega)^3$.
\begin{lemma} \label{restriction_operator}
There exists  a restriction operator $R_\ep$ acting from $H^1_0(\Omega)^3$ into $H^1_0(\Omega_\ep)^3$ such that
\begin{enumerate}
\item $v \in H^1_0(\Omega_\ep)^3\ \Rightarrow\ R_\ep v=v\hbox{  in }\Omega_\ep$ (elements of $H^1_0(\Omega_\ep)$ are continuated by $0$ to $\Omega$).
\item ${\rm div}(v)=0\hbox{  on }\Omega\ \Rightarrow\  {\rm div}(R_\ep v)=0\hbox{  in }\Omega_\ep$.
\item There exists a positive constant $C$, independent of $\ep$, such that
\begin{equation}\label{estim_restricted}
\begin{array}{l}
\|R_\ep v\|_{L^2(\Omega_\ep)^{3}}+ \ep\|D R_\ep v\|_{L^2(\Omega_\ep)^{3\times 3}} \leq C\left(\|v\|_{L^2(\Omega)}+\ep \|D v\|_{L^2(\Omega)^{3\times 3}}\right)\leq C \|v\|_{H^1_0(\Omega)^3}\,.
\end{array}
\end{equation}
\end{enumerate}
\end{lemma}
\begin{lemma}\label{lemma_ext_pressure} Let $q_\ep$ be a function in $L^2_0(\Omega_\ep)$. There exists a unique function $Q_\ep\in L^2_0(\Omega)$  which satisfies the following equality
\begin{equation}\label{ext_pressure_general}
\langle \nabla Q_\ep, v\rangle_{H^{-1},H^1_0(\Omega)}=\langle \nabla q_\ep, R_\ep v\rangle_{H^{-1},H^1_0(\Omega_\ep)},\quad\hbox{for every }v\in H^1_0(\Omega)^3.
\end{equation}
\end{lemma}
We denote by $P_\ep$ the extension of the pressure $p_\ep$ obtained by applying Lemma \ref{lemma_ext_pressure} and give the following result.
\begin{lemma}\label{estim_ext_pressure}Assume that the asymptotic regimes (\ref{RM}) and (\ref{body_forces}) and condition (\ref{condition_existence}) hold. Then there exists a positive constant $C$ independent of $\ep$, such that the following estimate    holds
\begin{equation}\label{estimates_p_ext}
\ep\|P_\ep\|_{L^2(\Omega)}+\ep\|\nabla P_\ep\|_{H^{-1}(\Omega)^3}\leq C.
\end{equation}
\end{lemma}
\begin{proof} From the definition (\ref{ext_pressure_general}) of the extension $P_\ep$ and  the variational formulation (\ref{var_form_v}), we get 
$$
\begin{array}{rl}  \displaystyle 
\langle \nabla P_\ep, v\rangle_{H^{-1}, H^1_0(\Omega)}=&\displaystyle -\int_{\Omega_\varepsilon}\rot(u_\ep)\cdot \rot(R_\ep v)\,dx+2N^2\int_{\Omega_\varepsilon}\rot(R_\ep v)\cdot w_\ep\,dx\\
\noame &\displaystyle-2\left({1\over \alpha}-N^2\right)\int_{\partial F_\varepsilon}(w_\ep\times n_\ep)\cdot R_\ep v\,d\sigma +\ep^{-1}\int_{\Omega_\ep}f\cdot R_\ep v\,dx.
\end{array}
$$
Applying Cauchy-Schwarz's inequality and taking into account  estimates of the velocity (\ref{estimates_u}),  microrotation (\ref{estimates_w}) and restricted operator  (\ref{estim_restricted}), we get 
$$\begin{array}{rl}\displaystyle\left|\int_{\Omega_\varepsilon}\rot(u_\ep)\cdot \rot(R_\ep v)\,dx\right|\leq&\displaystyle  C\|D u_\ep\|_{L^2(\Omega_\ep)^{3\times 3}}
\|D R_\ep v\|_{L^2(\Omega_\ep)^{3\times 3}}\leq {C\ep^{-1}}\|v\|_{H^1_0(\Omega)},\\
\noame
\displaystyle \left|\int_{\Omega_\varepsilon}\rot(R_\ep v)\cdot w_\ep\,dx\right|\leq &\displaystyle C\|D R_\ep v\|_{L^2(\Omega_\ep)^{3\times 3}} \|w_\ep\|_{L^2(\Omega_\ep)^{3}}\leq {C\ep^{-1}}\|v\|_{H^1_0(\Omega)},
\\
\noame
\displaystyle
\left|\ep^{-1}\int_{\Omega_\ep}f\cdot R_\ep v\,dx\right|\leq& \displaystyle \ep^{-1}\|f\|_{L^2(\Omega)^3}\|R_\ep v\|_{L^2(\Omega_\ep)^3}\leq {C \ep^{-1}}\|v\|_{H^1_0(\Omega)}.
\end{array}$$
From estimate (\ref{estimate_buo}) applied to $R_\ep v$ and using estimate (\ref{estim_restricted}), we deduce $\|R_\ep v\|_{L^2(\partial F_\ep)^{3}}\leq C\ep^{-{1\over 2}}$ and from estimate (\ref{Poincare}) with estimate (\ref{estimates_w}), we deduce   $\|w_\varepsilon\|_{L^2(\partial F_\ep)^{3}}\leq C\ep^{-{1\over 2}}$. Then, we get
$$\left|\int_{\partial F_\varepsilon}(w_\ep\times n_\ep)\cdot R_\ep v\,d\sigma\right|\leq C\|w_\ep\|_{L^2(\partial F_\ep)^{3}}\|R_\ep v\|_{L^2(\partial F_\ep)^{3}}\leq C\ep^{-{1\over 2}}\|R_\ep v\|_{L^2(\partial F_\ep)^{3}}\leq {C}\varepsilon^{-1}\|v\|_{H^1_0(\Omega)}. 
$$
This together with previous inequalities gives $|\langle \nabla P_\ep, v\rangle_{H^{-1}, H^1_0(\Omega)}|\leq {C \ep^{-1}}\|v\|_{H^1_0(\Omega)^3}$ and so $\|\nabla P_\ep\|_{H^{-1}(\Omega)^3}\leq C\ep^{-1}$.  By using the classical inequality  (see \cite{Temam})  
\begin{equation}\label{Necas}
\|P_\ep\|_{L^2(\Omega)}\leq C\|\nabla P_\ep\|_{H^{-1}(\Omega)^3},
\end{equation} we get estimate (\ref{estimates_p_ext}).
\end{proof}

\paragraph{A compactness result.}
Our aim is to describe the asymptotic behavior of the velocity $u_\ep$, microrotation $w_\ep$ and pressure $p_\ep$ of the fluid as $\ep$ tends  to $0$ taking into account the boundary of the obstacles. To do this we use the periodic unfolding method in perforated domains and the estimates given in the previous section. Thus, we briefly recall the definition of the unfolding operator and its main properties (for more details, see \cite{Ciora_donato_Griso, Dam} for fixed domains and \cite{Ciora_Don_Zak,Ciora_Don_Zak_0, Cioranescu_book}  for perforated domains).

In the sequel we will use the following notation:
\begin{itemize}
\item $\tilde \varphi$ for the zero extension outside $\Omega_\ep$ (resp. $\Omega)$ for any function $\varphi$ in $L^2(\Omega_\ep)$ (resp. $L^2(\Omega)$),
\item For $x\in \mathbb{R}^3$, we set $x=[x]_Y + \{x\}_Y$ where the integer part $[x]_Y$ belongs to the periodical net of $\mathbb{R}^3$ (i.e. the subgroup $\mathbb{Z}^3$) with respect to $Y$ and $\{x\}_Y=x-[x]_Y$ is the fractional part of $x$. Thus, for every $\ep>0$, the former decomposition implies that we also have $x=\ep\{x/\ep\}_Y+\ep[x/\ep]_Y$ for every $x\in \mathbb{R}^3$. 
\end{itemize}

\begin{definition}The unfolding operator $\mathcal{T}_\ep:L^2(\Omega_\ep)\to L^2(\mathbb{R}^3\times Y^*)$ is defined by 
$$
\mathcal{T}_\ep(\varphi)(x,y)=\tilde\varphi\left(\ep\left[{x\over \ep}\right]_Y+\ep y\right),\quad \forall\varphi \in L^2(\Omega_\ep),\quad \forall (x,y)\in\mathbb{R}^3\times Y^*.
$$
\end{definition}

\begin{proposition}\label{properties_UO}
The unfolding operator $\mathcal{T}_\ep$ has the following properties:
\begin{enumerate}
\item $\mathcal{T}_\ep$ is a linear operator.
\item $\mathcal{T}_\ep(\varphi\phi)=\mathcal{T}_\ep(\varphi)\mathcal{T}_\ep(\phi),\ \forall\,\varphi,\phi\in L^2(\Omega_\ep)$.
\item $\mathcal{T}_\ep(\varphi_\ep)(x,y)=\varphi(y),\ \forall\,(x,y)\in \mathbb{R}^3\times Y^*,\ \forall\,\varphi\in L^2(Y^*)$ a $Y$-periodic function with $\varphi_\ep(x)=\varphi\left({x\over\ep}\right)$.
\item $\|\mathcal{T}_\ep(\varphi)\|_{L^2(\mathbb{R}^3\times Y^*)}=|Y|^{1\over 2}\|\varphi\|_{L^2(\Omega_\ep)},\  \forall\varphi\in L^2(\Omega_\ep)$.
\item $\nabla_y \mathcal{T}_\ep(\varphi)(x,y)=\ep \mathcal{T}_\ep(\nabla_x \varphi)(x,y),\ \forall\,(x,y)\in \mathbb{R}^3\times Y^*, \ \forall\,\varphi\in H^1(\Omega_\ep)$.
\item $\mathcal{T}_\ep(\varphi)\in L^2(\mathbb{R}^3; H^1(Y^*)),\  \forall\,\varphi\in H^1(\Omega_\ep)$.
\item Let $\varphi_\ep$ be in $L^2(\Omega)$ such that $\tilde \varphi_\ep\to \varphi$   in $L^2(\Omega)$. Then $\mathcal{T}_\ep(\varphi_\ep)\to \tilde\varphi$ in $L^2(\mathbb{R}^3\times Y^*)$.
\end{enumerate}
\end{proposition}

\begin{proposition}\label{Prop_cong_uhat}Let $\varphi_\ep$ be a sequence such that 
$$\|\varphi_\ep\|_{L^2(\Omega_\ep)^3}+\ep\|D \varphi_\ep\|_{L^2(\Omega_\ep)^{3\times 3}}\leq C\,.$$
Then, there exists $\hat \varphi\in L^2(\Omega;H^1_{\rm per}(Y^*)^3)$ such that 
$$
\mathcal{T}_\ep(\varphi_\ep)\rightharpoonup \hat \varphi\quad\hbox{ in }L^2(\Omega;H^1(Y^*)^3),\quad \ep\mathcal{T}_\ep(D \varphi_\ep)\rightharpoonup D_y\hat \varphi\quad\hbox{ in }L^2(\Omega\times Y^*)^{3\times 3}.
$$
\end{proposition}
In a similar way,  it is introduced in \cite{Ciora_Don_Zak, Cioranescu_book} the unfolding operator
on the boundary of the holes $\partial F_\ep$. 

\begin{definition}The unfolding boundary operator $\mathcal{T}_\ep^b(\varphi)\in L^2(\mathbb{R}^3\times \partial F)$ is defined by 
\begin{equation}\label{boundary_unfolding_operator}
\mathcal{T}_\ep^b(\varphi)(x,y)=\tilde\varphi\left(\ep\left[{x\over \ep}\right]_Y+\ep y\right),\quad \forall\,\varphi\in L^2(\partial F_\ep),\quad \forall (x,y)\in\mathbb{R}^3\times \partial F.
\end{equation}
\end{definition}
We remark that if $\varphi\in H^1(\Omega_\ep)$ and $\varphi=0$ on $\partial \Omega$, one has $\mathcal{T}_\ep^b(\varphi)=\mathcal{T}_\ep(\varphi)$ on $\partial F$. 
The next results reformulate the properties given above in the case of functions
defined on the boundary of the holes $\partial F_\ep$. 
\begin{proposition}\label{properties_BUO}
The unfolding boundary operator $\mathcal{T}^b_\ep$ has the following properties:
\begin{enumerate}
\item $\mathcal{T}_\ep^b$ is linear.
\item $\mathcal{T}_\ep^b(\varphi\phi)=\mathcal{T}_\ep^b(\varphi)\mathcal{T}_\ep^b(\phi),\ \forall\,\varphi,\phi\in L^2(\partial F_\ep)$.
\item $\mathcal{T}_\ep^b(\varphi_\ep)(x,y)=\varphi(y),\quad \forall\,(x,y)\in \mathbb{R}^3\times \partial F,\ \forall\,\varphi\in L^2(\partial F)$ a $Y$-periodic function with $\varphi_\ep(x)=\varphi\left({x\over\ep}\right)$.
\item $\|\mathcal{T}_\ep^b(\varphi)\|_{L^p(\mathbb{R}^3\times \partial F)}=(\ep|Y|)^{1\over 2}\|\varphi\|_{L^p(\partial F_\ep)},\  \forall\varphi\in L^2(\partial F_\ep)$.
\item $\lim_{\ep\to 0}\int_{\mathbb{R}^3\times \partial F}\mathcal{T}_\ep^b(\varphi)(x,y)dxd\sigma(y)=|\partial F|\int_\Omega \varphi(x)\,dx,\ \forall\,\varphi\in H^1(\Omega)$.
\item $\mathcal{T}_\ep^b(\varphi)\to \tilde\varphi$ strongly in $L^2(\mathbb{R}^3\times\partial F),\ \forall\,\varphi\in H^1_0(\Omega)$.
\end{enumerate}
\end{proposition}

Next, we give some compactness results about the behavior
of  the extended functions $(U_\ep, W_\ep,  P_\ep)$ and the unfolding functions $(\mathcal{T}_\ep(u_\ep),  \mathcal{T}_\ep(w_\ep), \mathcal{T}_\ep( P_\ep))$ by assuming the  {\it a priori} estimates given in Lemmas \ref{Estimates_extended_lemma}, \ref{lemma_ext_vel} and \ref{estim_ext_pressure}.
\begin{proposition} \label{Prop_convergences_uwp}Assume that the asymptotic regimes (\ref{RM}) and (\ref{body_forces}) and condition (\ref{condition_existence}) hold. Then,  for a subsequence of $\ep$ still denoted by $\ep$, we have that 
\begin{enumerate}
\item (Velocity) there exist $u\in L^2(\Omega)^3$ and $\hat u\in L^2(\Omega; H^1_{\rm per}(Y^*)^3)$ such that 
\begin{equation}\label{conv_hat_u_0}
\begin{array}{l}
\ep^{-1}U_\ep\rightharpoonup u\hbox{ in }L^2(\Omega)^3,
\end{array}
\end{equation}
\begin{equation}\label{conv_hat_u_1}
\begin{array}{l}
\ep^{-1}\mathcal{T}_\ep(u_\ep)\rightharpoonup \hat u\hbox{ in }L^2(\Omega;H^1(Y^*))^3,\quad \mathcal{T}_\ep(Du_\ep)\rightharpoonup D_y\hat u\hbox{ in }L^2(\Omega\times Y^*)^{3\times 3},
\end{array}
\end{equation}
\begin{equation}\label{conv_hat_u_1_b}
\ep^{-1}\mathcal{T}^b_\ep(u_\ep)\rightharpoonup \hat u\hbox{ in }L^2(\Omega;H^{1\over 2}(\partial F))^3,
\end{equation}
\begin{equation}\label{conv_hat_u_1_rot}
\begin{array}{l}
 \mathcal{T}_\ep(\rot(u_\ep))\rightharpoonup \rot_y(\hat u)\hbox{ in }L^2(\Omega\times Y^*)^{3},\end{array}
\end{equation}
and moreover,  the following conditions hold
\begin{eqnarray}
u(x)=\int_{Y^*}\hat u(x,y)\,dy&&\hbox{in }\Omega,\label{relation_u_hatu}\\
\noame
\hat u(x,y)\cdot n(y)=0&&\hbox{on }\Omega\times \partial F,\label{slip_hat_u}\\
\noame
{\rm div}_y\hat u(x,y)=0&&\hbox{in }\Omega\times Y^*,\label{div_hat_u}\\
\noame
{\rm div}_x\left(\int_{Y^*}\hat u(x,y)\,dy\right)=0&&\hbox{in }\Omega,\label{divxhatu}\\
\noame
\left(\int_{Y^*}\hat u(x,y)\,dy\right)\cdot n(y)=0&&\hbox{in }\Omega,\label{divxhatunormal}
\end{eqnarray}

\item (Microrotation) there exist $w\in L^2(\Omega)^3$ and $\hat w\in L^2(\Omega; H^1_{\rm per}(Y^*)^3)$ such that 
\begin{equation}\label{conv_hat_w_0}
\begin{array}{l}
W_\ep\rightharpoonup w\hbox{ in }L^2(\Omega)^3,\end{array}
\end{equation}
\begin{equation}\label{conv_hat_w_1}
\begin{array}{l}
\mathcal{T}_\ep(w_\ep)\rightharpoonup \hat w\hbox{ in }L^2(\Omega;H^1(Y^*))^3,\quad\ep\mathcal{T}_\ep(Dw_\ep)\rightharpoonup D_y\hat w\hbox{ in }L^2(\Omega\times Y^*)^{3\times 3},
\end{array}
\end{equation}
\begin{equation}\label{conv_hat_w_1_b}
\mathcal{T}^b_\ep(w_\ep)\rightharpoonup \hat w\hbox{ in }L^2(\Omega;H^{1\over 2}(\partial F))^3,
\end{equation}
\begin{equation}\label{conv_hat_w_1_rot}
\ep\mathcal{T}_\ep(\rot(w_\ep))\rightharpoonup \rot_y(\hat u)\hbox{ in }L^2(\Omega\times Y^*)^{3},
\end{equation}
and moreover,  the following conditions hold
\begin{eqnarray}
w(x)=\int_{Y^*}\hat w(x,y)\,dy&&\hbox{in }\Omega,\label{relation_w_hatu}\\
\hat w(x,y)\cdot n(y)=0&&\hbox{on }\Omega\times \partial F,\label{slip_hat_w}
\end{eqnarray}
\item (Pressure) there exist $p\in L^2_0(\Omega)^3$ such that 
\begin{equation}\label{conv_p}
\ep P_\ep\to p\hbox{ in }L^2(\Omega),\quad \ep \mathcal{T}_\ep(P_\ep)\to p\hbox{ in }L^2(\Omega\times Y^*).\quad 
\end{equation}
\end{enumerate}
\end{proposition}
\begin{proof} We start proving {\it 1}.  From estimates for the extended velocity (\ref{estimates_u_ext}) we deduce convergence (\ref{conv_hat_u_0}). Taking into account the a priori estimates for the velocity (\ref{estimates_u})  and using Proposition  \ref{Prop_cong_uhat}, we deduce  convergences given in  the (\ref{conv_hat_u_1}). Convergence (\ref{conv_hat_u_1_b}) is
straightforward from the definition (\ref{boundary_unfolding_operator}) and the Sobolev injections. Finally, taking into account Proposition \ref{properties_UO}$_{1,5}$, we have
$$\mathcal{T}_\ep\left(\partial_{x_i}u_{\ep,j}-\partial_{x_j}u_{\ep,i}\right)= \mathcal{T}_\ep\left(\partial_{x_i}u_{\ep,j}\right)-\mathcal{T}_\ep\left(\partial_{x_j}u_{\ep,i}\right)=\ep^{-1}(\partial_{y_i}\mathcal{T}_\ep(u_{\ep,j})-\partial_{y_j}\mathcal{T}_\ep(u_{\ep,i})),\quad \forall\, i,j=1,2,3,\, i<j,$$
and so $\mathcal{T}_\ep(\rot(u_\ep))=\ep^{-1}\rot_y\left(\mathcal{T}_\ep(u_\ep)\right)$, which from convergence (\ref{conv_hat_u_1}) implies  (\ref{conv_hat_u_1_rot}).

In order to prove the boundary conditions (\ref{slip_hat_u}), let us take $\varphi\in \mathcal{D}(\Omega)$ and from $u_\ep\cdot n_\ep=0$ on $\partial F_\ep$, we have
$$\int_{\partial F_\ep}(u_\ep\cdot n_\ep)\varphi\,d\sigma(x)=0.$$
By applying the unfolding boundary and using Proposition \ref{properties_BUO}$_{2,3,4}$, we get
$$0=\int_{\partial F_\ep}\!(u_\ep\cdot n_\ep)\varphi\,d\sigma(x)=\ep^{-1}\!\int_{\mathbb{R}^3\times \partial F}\!\left(\mathcal{T}_\ep^b(u_\ep)\cdot  n\right)\mathcal{T}_\ep^b(\varphi)\,dx\,d\sigma(y).$$
Passing to the limit when $\ep$ tends to zero, from convergence (\ref{conv_hat_u_1_b})  and Proposition \ref{properties_BUO}$_{6}$, we obtain
$$0=\int_{\Omega\times \partial F}(\hat u(x,y)\cdot n(y))\varphi\,dxd\sigma(y)=\int_{\Omega}\left(\int_{\partial F}\hat u(x,y)\cdot n(y)\,dy\right)\varphi(x)\,dx\,,$$
which implies (\ref{slip_hat_u}).

In order to prove relation (\ref{div_hat_u}), let us observe that ${\rm div}(u_\ep)=0$ implies $\mathcal{T}_\ep({\rm div}(u_\ep))=0$. But from Proposition \ref{properties_UO}$_{1,5}$, we have
$$\mathcal{T}_\ep({\rm div}(u_\ep))=\sum_{i=1}^3\mathcal{T}_\ep\left(\partial_{x_i}u_{\ep,i}\right)= \ep^{-1}\sum_{i=1}^3\partial_{y_i}\mathcal{T}_\ep\left(u_{\ep,i}\right)=\ep^{-1}{\rm div}_y(\mathcal{T}_\ep(u_\ep))$$
and so $\ep^{-1}{\rm div}_y(\mathcal{T}_\ep(u_\ep))=0$. Passing to the limit as $\ep$ tends to zero in the last equality we get (\ref{div_hat_u}).

In order to prove (\ref{slip_hat_u}) and (\ref{divxhatu}), multiplying ${\rm div}(u_\ep)=0$ by $\ep^{-1}\varphi$ with  $\varphi$ in $\mathcal{D}(\Omega)$ and  using $u_\ep\cdot n_\ep=0$ on $\partial F_\ep$, we have 
$$0=\int_{\Omega_\ep} \ep^{-1}{\rm div}(u_\ep)\varphi\,dx=\int_{\Omega_\ep} \ep^{-1}u_\ep\cdot \nabla \varphi\,dx.$$
By applying the unfolding, we get 
$$\int_{\mathbb{R}^3\times Y^*} \ep^{-1}\mathcal{T}_\ep(u_\ep)\cdot
\mathcal{T}_\ep(\nabla \varphi)\,dxdy=0.$$
We pass to the limit as $\ep$ tends to zero and we get
$$\int_{\Omega\times Y^*} \hat u(x,y)\cdot
\nabla \varphi(x)\,dxdy=0,$$
and so 
$$\int_{\Omega}{\rm div}_x\left(\int_{Y^*} \hat u(x,y)\,dy\right)
 \varphi(x)\,dx=0,\quad \forall\,\varphi\in \mathcal{D}(\Omega),$$
 which implies (\ref{divxhatu}) and (\ref{divxhatunormal}). 
 
 Finally, we prove (\ref{relation_u_hatu}). From   Proposition \ref{properties_UO} and  taking $\varphi_\ep=\ep^{-1}\varphi$ with $\varphi$ in $\mathcal{D}(\Omega)$, we have 
 $$\ep^{-1}\int_{\Omega_\ep}u_\ep(x)\cdot \varphi(x)\,dx={1\over \ep |Y|}\int_{\mathbb{R}^3\times Y^*}\mathcal{T}_\ep(u_\ep)(x,y)\cdot \mathcal{T}_\ep(\varphi)(x,y)\,dxdy.$$
 By using the extension of the velocity, we have 
  $$\ep^{-1}\int_{\Omega}U_\ep(x)\cdot \varphi(x)\,dx=\ep^{-1}\int_{\mathbb{R}^3\times Y^*}\mathcal{T}_\ep(u_\ep)(x,y)\cdot \mathcal{T}_\ep(\varphi)(x,y)\,dxdy,$$
  and passing to the limit by using convergences (\ref{conv_hat_u_0}) and (\ref{conv_hat_u_1}), we get
$$\int_{\Omega}u(x)\cdot \varphi(x)\,dx=\int_{\Omega}\left(\int_{Y^*}\hat u(x,y)\,dy\right)\varphi(x)\,dx,$$
 which implies   property (\ref{relation_u_hatu}).

The proof of {\it 2} is similar, so we omit it. 

We finish with the proof of {\it 3}. Thus, estimates given in Lemma \ref{estim_ext_pressure} imply, up to a subsequence, the existence of $p\in L^2_0(\Omega)$ such that 
\begin{equation}\label{weak_conv_p}\ep P_\ep\rightharpoonup p\hbox{ in }L^2(\Omega),\quad \ep \nabla P_\ep\rightharpoonup p\hbox{ in }H^{-1}(\Omega)^3.
\end{equation} Moreover, following \cite{Tartar} it can be proved that this convergence is in fact strong.  To prove this, let $\sigma_\ep\in H^1_0(\Omega)^3$ such that 
\begin{equation}\label{conv_sigma}
\sigma_\ep\rightharpoonup \sigma \hbox{ in }H^1_0(\Omega)^3.
\end{equation}
Then, we have 
$$\left|\langle\ep\nabla P_\ep,\sigma_\ep\rangle_{H^{-1},H^1_0(\Omega)^3}\right|\leq \left|\langle \ep\nabla P_\ep,\sigma_\ep-\sigma\rangle_{H^{-1},H^1_0(\Omega)^3}\right|+\left|\langle\ep\left(\nabla P_\ep-\nabla p\right), \sigma\rangle_{H^{-1},H^1_0(\Omega)^3}\right|.$$
On the one hand, using first convergence  in (\ref{weak_conv_p}), we have
$$\left|\langle\ep\left(\nabla P_\ep-\nabla p\right),\sigma\rangle_{H^{-1},H^1_0(\Omega)^3}\right|=\left|\int_\Omega\ep\left(P_\ep-p\right){\rm div}(\sigma)\,dx\right|\to 0,\quad\hbox{as }\ep\to 0.$$
On the other hand,  from estimate (\ref{estim_restricted}) and (\ref{estimates_p_ext}), we have that
$$\begin{array}{rl}
\displaystyle\left|\langle\ep \nabla P_\ep,\sigma_\ep-\sigma\rangle_{H^{-1},H^1_0(\Omega)^3}\right|=& \displaystyle \left|\langle \ep \nabla P_\ep, R_\ep (\sigma_\ep-\sigma)\rangle_{H^{-1},H^1_0(\Omega_\ep)^3}\right|\\
\noame
\leq &\displaystyle C\left(\| \sigma_\ep-\sigma\|_{L^2(\Omega)^3}+\ep\|D(\sigma_\ep-\sigma)\|_{L^2(\Omega)^{3}}\right)\to 0,\hbox{ as }\ep\to 0,
\end{array}$$
by virtue of (\ref{conv_sigma}) and the Rellich theorem. This implies that $\nabla P_\ep\to \nabla p$ in $H^{-1}(\Omega)^3$, which together with  inequality (\ref{Necas}), implies the strong convergence of the pressure $P_\ep$ given in (\ref{conv_p}).
This convergence and Proposition \ref{properties_UO}$_7$ imply the strong convergence of $\ep\mathcal{T}_\ep(P_\ep)$ to $p$ in $L^2(\Omega\times Y^*)$.
\end{proof}
\paragraph{Obtaining the limit system.} We use the results of the previous sections to prove Theorem \ref{thm_darcy_law} describing the asymptotic behavior of the solution of the micropolar system  (\ref{system_1_1})-(\ref{BC_exterior}). To do this, we first give the existence and uniqueness result for micropolar local problem (\ref{local_problems}).

\begin{lemma}\label{lema_cell_problems}Assume that condition (\ref{condition_existence}) holds. Then,  for every $k=1,2$ and $i=1,2,3$, there exists 
a unique solution  $(u^{i,k}, w^{i,k}, \pi^{i,k})\in V^0_Y\times V_Y\times L^2_0(Y^*)$ of  the local problem (\ref{local_problems}).
 \end{lemma}
 \begin{proof} Similarly to Proposition \ref{prop_form_var}, sufficiently regular solutions of (\ref{local_problems}) satisfy the following weak formulation: For $i=1,2,3$, $k=1,2$ find $(u^{i,k},w^{i,k},\pi^{i,k})\in V_Y^0\times V_Y\times L^2_0(Y^*)$   such that
$$
\begin{array}{l}\displaystyle
\int_{Y^*}\rot(u^{i,k})\cdot \rot(\varphi)\,dx-\int_{Y^*}\pi^{i,k} {\rm div}(\varphi)\,dx-2N^2\int_{Y^*}\rot(\varphi)\cdot w^{i,k}\,dx
\\
\noame
\displaystyle+2\left({1\over \alpha}-N^2\right)\int_{\partial F}(w^{i,k}\times n)\cdot\varphi\,d\sigma= \int_{Y^*}e_i\delta_{1k}\cdot \varphi\,dx\,,\  \forall \varphi\in V_Y\,,
\end{array}
$$
$$
\begin{array}{l}\displaystyle
R_c\int_{Y^*}\rot(w^{i,k})\cdot \rot(\psi)\,dx+R_c\int_{Y^*}{\rm div}(w^{i,k})\cdot {\rm div}(\psi)\,dx+4N^2\int_{Y^*}w^{i,k}\cdot\psi\,dx\\
\noame\displaystyle
\qquad-2N^2\int_{Y^*}\rot(\psi)\cdot u^{i,k}\,dx + 2N^2(\beta-1)\int_{\partial F}(u^{i,k}\times n)\cdot \psi\,d\sigma
= \int_{Y^*}e_i\delta_{2k}\cdot \psi\,dx\,,\quad \forall \psi\in V_Y\,.
\end{array}
$$
 Thus, following the lines of the proof of Theorem \ref{lema_coercive},  we can introduce the following mixed variational form:

For $i=1,2,3$, $k=1,2$ find $(u^{i,k},w^{i,k},\pi^{i,k})\in V_Y^0\times V_Y\times L^2_{0,{\rm per}}(Y^*)$ such that
$$
\begin{array}{rcll}
\mathcal{A}_Y(u^{i,k}, w^{i,k};\varphi,\psi)+\mathcal{B}_Y((\varphi,\psi),\pi^{i,k})&=&\mathcal{L}^{i,k}_Y(\varphi,\psi)  & \forall\, (\varphi,\psi)\in V_Y\times V,
\\
\noame
\mathcal{B}_\ep((u^{i,k}, w^{i,k}),q^{i,k})&=&0&    \forall\, q^{i,k}\in L^2_{0,{\rm per}}(Y^*)\,,
\end{array}
$$
where 
\begin{equation}\label{Def_A_cell}\begin{array}{rl}\displaystyle
\mathcal{A}_Y(u^{i,k},w^{i,k};\varphi,\psi)=&\displaystyle 
\int_{Y^*}\rot_y(u^{i,k})\cdot \rot_y(\varphi)\,dx-2N^2\int_{Y^*}\rot_y(\varphi)\cdot w^{i,k}\,dy\\
\noame&
\displaystyle +R_c\int_{Y^*}\rot_y(w^{i,k})\cdot \rot_y(\psi)\,dy+R_c\int_{Y^*}{\rm div}_y(w^{i,k})\cdot {\rm div}_y(\psi)\,dy
\\
\noame& \displaystyle
-2N^2\int_{Y^*}\rot_y(\psi)\cdot u^{i,k}\,dy+4N^2\int_{Y^*}w^{i,k}\cdot\psi\,dy\\
\noame &\displaystyle
+2\left({1\over \alpha}-N^2\right)\int_{\partial F}(w^{i,k}\times n)\cdot\varphi\,d\sigma(y)+ 2N^2(\beta-1)\int_{\partial F}(u^{i,k}\times n)\cdot \psi\,d\sigma(y)\,,
\end{array}
\end{equation}
 and 
$$
\begin{array}{rl}\displaystyle
\mathcal{B}_Y((\varphi,\psi),\pi^{i,k})=-\int_{Y^*}\pi^{i,k}{\rm div}_y \varphi\,dy,\quad \mathcal{L}^{i,k}_Y(\varphi,\psi)=&\displaystyle \int_{Y^*}e_i\delta_{1k}\cdot \varphi\,dy+
\int_{Y^*}e_i\delta_{2k}\cdot \psi\,dy\,.
\end{array}
$$
We denote
$$\|(\varphi,\psi)\|_{V_Y^0\times V_Y}=(\|D \varphi\|^2_{L^2(Y^*)^{3\times 3}}+\|D \psi\|^2_{L^2(Y^*)^{3\times 3}})^{1\over 2}.$$ 
Following the proof of Theorem \ref{lema_coercive} and taking into account that in $Y^*$ the trace inequality (\ref{estimate_buo}) holds with constant $C_t$ instead of $\ep^{-1} C_t$,  the Poincar\'e inequality (\ref{Poincare}) with constant $C_p$ instead of $\ep C_p$, the trace inequality (\ref{Unfolding_boundary_1}) with constant  $C_{pt}^{1\over 2}$ instead of  $(\ep C_{pt})^{1\over 2}$ and the Gaffney inequality (\ref{Maxwell_inequality}) with constant $C_g$, 
it is not difficult to prove that the bilinear form $\mathcal{A}_Y$, $\mathcal{B}_Y$ and $\mathcal{L}_Y$  are continuous bilinear forms on $(V_Y^0\times V_Y)^2$, $V_Y\times V_Y\times L^2_{0,{\rm per}}(Y^*)$  and $V_Y\times V_Y$ respectively. . Moreover,  under condition  (\ref{condition_existence}) it follows from the proof of the coercivity of $\mathcal{A}_\ep$ in Theorem \ref{lema_coercive} that $\mathcal{A}_Y$ satisfies 
$$
\mathcal{A}_Y(\varphi,\psi;\varphi,\psi)\geq A \|D \varphi\|^2_{L^2(\Omega_\ep)^{3\times 3}}+B \|D \psi\|^2_{L^2(\Omega_\ep)^{3\times 3}}
\geq \min\{A,B\}\|(\varphi,\psi)\|^2_{V_Y^0\times V_Y},
$$
with positive constants $A$, $B$ given by (\ref{AB}). Moreover, it can be proved  the inf-sup condition 
$$\exists\, \delta>0,\quad\hbox{such that}\quad \sup_{(\varphi,\psi)\in V_Y\times V_Y}{\int_{Y^*}{\rm div}\,\varphi\,q^{i,k}\,dx\over (\|\varphi\|^2_{V_Y}+\|\varphi\|^2_{V_Y})^{1\over 2}}
\geq
\delta,$$
which ends the proof.
\end{proof}

\begin{proof}[Proof of Theorem \ref{thm_darcy_law}] We divide the proof in two steps.

{\it Step 1}. In this step we prove that the whole sequences $(\ep^{-1}\mathcal{T}_\ep(u_\ep), \mathcal{T}_\ep(w_\ep))$  and $\ep \mathcal{T}_\ep(P_\ep)$ converge weakly in $L^2(\Omega;H^1(Y^*)^3)\times L^2(\Omega;H^1(Y^*)^3)$ to $(\hat u, \hat w)$   and strongly to $p$ in $L^2(\Omega\times Y^*)$ respectively, where the triplet   $(\hat u,\hat w, p)\in L^2(\Omega;H^1_{\rm per}(Y^*)^3)\times L^2(\Omega;H^1_{\rm per}(Y^*)^3)\times L^2_0(\Omega)$  is the unique solution of the following homogenized system
\begin{eqnarray}
-\Delta_y \hat u+ \nabla_y\hat q-2N^2{\rm rot}_y(\hat w)= f -\nabla_x p &\quad\hbox{in}\quad\Omega\times Y^*,\label{system_1_1_hat}\\
\noame
-R_c\Delta_y \hat w+4N^2\hat w-2N^2{\rm rot}_y(\hat u) = g&\quad\hbox{in}\quad\Omega\times Y^*,\label{system_1_3_hat}\\
\noame
{\rm div}_y(\hat u)=0&\quad\hbox{in}\quad\Omega\times Y^*,\label{system_1_2_hat}\\
{\rm div}_x\left(\int_{Y^*}\hat u\,dy\right)=0&\hbox{in }\Omega,\label{divxhatu_system}\\
\left(\int_{Y^*}\hat u\,dy\right)\cdot n=0&\hbox{in }\Omega,\label{divxhatunormal_system}
\end{eqnarray}
 where $\hat q\in  L^2(\Omega;L^2_{0,{\rm per}}(Y^*)^3)$, and the boundary conditions
\begin{eqnarray}
\displaystyle \hat w \times n={\alpha\over 2} \rot_y(\hat u)\times n & \hbox{ on }\Omega\times \partial F,\label{BC_holes_1_hat}\\
\noame
\rot_y(\hat w)\times n =\displaystyle{2N^2\over R_c} \beta (\hat u\times n)& \hbox{ on }\Omega\times \partial F,\label{BC_holes_2_hat}\\
\noame
\hat u\cdot n =0 & \hbox{ on }\Omega\times \partial F,\label{BC_holes_3_hat}\\
\noame
\hat w\cdot n =0 & \hbox{ on }\Omega\times \partial F.\label{BC_holes_4_hat}
\end{eqnarray}
By taking into account Proposition \ref{Prop_convergences_uwp}, we have that (\ref{system_1_2_hat})-(\ref{divxhatunormal_system}) and (\ref{BC_holes_3_hat})-(\ref{BC_holes_4_hat})   hold. Below, we prove the rest of them.

First, we prove (\ref{system_1_1_hat}). To do this, we first  take as test function in (\ref{var_form_v}) the following function $\varphi_\ep(x)=\ep \phi(x) \Phi(x/\ep)$, where $\phi\in  \mathcal{D}(\Omega)$ and $\Phi \in H^1_{\rm per}(Y^*)^3$ with $\Phi\cdot n=0$ on $\partial F$ and ${\rm div}_y\Phi=0$ in $Y^*$. Then, we have
$$\begin{array}{l}\displaystyle
\int_{\Omega_\varepsilon} \rot( u_\ep)\cdot \rot(\varphi_\ep)\,dx-\int_{\Omega_\ep}p_\ep\,{\rm div}(\varphi_\ep)\,dx
-2N^2 \int_{\Omega_\varepsilon}\rot(\varphi_\ep)\cdot  w_\ep\,dx\\
\noame
\displaystyle+2 \left({1\over \alpha}-N^2\right)\int_{\partial F_\varepsilon}(w_\ep\times n_\ep)\cdot\varphi_\ep\,d\sigma(x)
= \ep^{-1}\int_{\Omega_\ep}f\cdot \varphi_\ep\,dx\,.
\end{array}$$
Let us observe   
\begin{equation}\label{div_rot}\rot(\varphi_\ep)=\ep \nabla\phi\times\Phi\left({\cdot\over \ep}\right)+ \phi\, \rot_y\left(\Phi\left({\cdot\over \ep}\right)\right),\quad {\rm div}(\varphi_\ep)=\ep \nabla\phi\cdot\Phi\left({\cdot\over \ep}\right)+ \phi\, {\rm div}_y\left(\Phi\left({\cdot\over \ep}\right)\right).
\end{equation}
Hence, by Proposition \ref{properties_UO}$_{2,3,7}$, we have   $\mathcal{T}_\ep(\varphi_\ep)=\ep \mathcal{T}_\ep(\phi) \Phi$ and convergences
\begin{equation}\label{property_phi}
\begin{array}{c}\displaystyle \mathcal{T}_\ep(\phi) \Phi\to  \phi \Phi\ \hbox{in }L^2(\Omega\times Y^*)^3,\  \mathcal{T}_\ep(\varphi_\ep)\to 0\ \hbox{in }L^2(\Omega\times Y^*)^3,\ \mathcal{T}_\ep(\rot(\varphi_\ep))
\rightharpoonup \phi\, \rot_y(\Phi)\ \hbox{in }L^2(\Omega\times Y^*)^3.
\end{array}
\end{equation}

By applying the unfolding to the variational formulation, taking into account the extension of the pressure and using Propositions \ref{properties_UO}$_{2,3}$ and \ref{properties_BUO}$_{2,3}$, we get
$$\begin{array}{l}\displaystyle
\int_{\mathbb{R}^3 \times Y^*}  \mathcal{T}_\ep(\rot(u_\ep))\cdot \mathcal{T}_\ep(\rot(\varphi_\ep))\,dxdy-\int_{\mathbb{R}^3\times Y^*}\mathcal{T}_\ep(P_\ep)\,\mathcal{T}_\ep({\rm div}(\varphi_\ep))\,dxdy
-2N^2\int_{\mathbb{R}^3 \times Y^*} \mathcal{T}_\ep(\rot(\varphi_\ep))\cdot \mathcal{T}_\ep(w_\ep) \,dxdy
\\
\noame\displaystyle
+2 \left({1\over \alpha}-N^2\right)\int_{\mathbb{R}^3\times \partial F}\left(\mathcal{T}^b_\ep(w_\ep)\times n\right)\cdot(\mathcal{T}^b_\ep(\phi) \Phi)\,dxd\sigma(y)= \int_{\mathbb{R}^3\times Y^*}\mathcal{T}_\ep(f)\cdot (\mathcal{T}_\ep(\phi) \Phi)\,dxdy\,.
\end{array}$$
Next, we pass to the limit in every terms of the previous variational formulation:

\begin{itemize}
\item First term. From convergence (\ref{conv_hat_u_1_rot}) and (\ref{property_phi}), we have 
$$\begin{array}{rl}
\displaystyle \int_{\mathbb{R}^3 \times Y^*}  \mathcal{T}_\ep(\rot(u_\ep))\cdot \mathcal{T}_\ep(\rot(\varphi_\ep))\,dxdy &\to\displaystyle \int_{\Omega \times Y^*}\rot_y(\hat u(x,y))\cdot (\phi(x)\rot_y(\Phi(y)))\,dxdy\\
\noame & \displaystyle =\int_{\Omega \times Y^*}\rot_y(\hat u(x,y))\cdot\rot_y(\phi(x)\Phi(y))\,dxdy.
\end{array}$$
\item Second term. We use (\ref{div_rot}) and Proposition \ref{properties_UO}$_{1,2,7}$, the fact that ${\rm div}_y(\Phi)=0$ in $Y^*$ and convergence (\ref{conv_p}), 
$$\begin{array}{rl}
\displaystyle \int_{\mathbb{R}^3\times Y^*}\mathcal{T}_\ep(P_\ep)\mathcal{T}_\ep({\rm div}(\varphi_\ep))\,dxdy& =
\displaystyle\int_{\mathbb{R}^3\times Y^*}\ep \mathcal{T}_\ep(P_\ep)\,\mathcal{T}_\ep(\nabla \phi\cdot \Phi)\,dxdy\\
\noame
 &\to\displaystyle \int_{\Omega\times Y^*}p(x)\nabla\phi(x)\cdot\Phi(y)\,dx=\int_{\Omega\times Y^*}p(x){\rm div}_x(\phi(x) \Phi(y))\,dxdy\\
\noame
&\displaystyle =-\int_{\Omega\times Y^*}\nabla_xp(x)(\phi(x) \Phi(y))\,dxdy\,.
\end{array}
$$
\item Third term. From convergences (\ref{conv_hat_w_1}) and (\ref{property_phi}), we get 
$$\begin{array}{rl}\displaystyle
\int_{\mathbb{R}^3\times Y^*} \mathcal{T}_\ep(\rot(\varphi_\ep))\cdot \mathcal{T}_\ep(w_\ep) \,dxdy\to &\displaystyle
\int_{\Omega \times Y^*}  \phi(x) \rot_y(\Phi(y))\cdot \hat w(x,y)\,dxdy\\
\noame=&\displaystyle
\int_{\Omega \times Y^*}  \rot_y(\phi(x) \Phi(y))\cdot \hat w(x,y)\,dxdy.
\end{array}$$
\item Fourth term. We use convergence (\ref{conv_hat_w_1_b})  and (\ref{property_phi}),
$$\int_{\mathbb{R}^3\times \partial F}\left(\mathcal{T}^b_\ep(w_\ep)\times n\right)\cdot(\mathcal{T}^b_\ep(\phi) \Phi)\,dxd\sigma(y)\to \int_{\Omega\times \partial F}\left(\hat w(x,y)\times n(y)\right)\cdot(\phi(x) \Phi(y))\,dxd\sigma(y).$$

\item Fifth term. From Proposition \ref{properties_UO}$_{7}$ and convergence (\ref{property_phi}),  we have 
$$\int_{\mathbb{R}^3\times Y^*}\mathcal{T}_\ep(f)\cdot (\mathcal{T}_\ep(\phi) \Phi)\,dxdy\to \int_{\Omega\times Y^*}f(x)\cdot (\phi(x)\Phi(y)) \,dxdy.$$
\end{itemize}
Therefore, taking into account the previous convergences and denoting  $\varphi(x,y)=\phi(x)\Phi(y)$, we obtain
$$\begin{array}{l}\displaystyle
\!\!\!\!\int_{\Omega \times Y^*}\!\!\!\!\rot_y(\hat u(x,y))\cdot \rot_y(\varphi(x,y))dxdy+\int_{\Omega}\nabla_x p(x)\, \varphi(x,y)dxdy
-2N^2\int_{\Omega \times Y^*}\!\!\!\!  \rot_{y}(\varphi(x,y))\cdot \hat w(x,y)dxdy
\\
\noame\displaystyle
\qquad
+2 \left({1\over \alpha}-N^2\right)\int_{\Omega\times \partial F}(\hat w(x,y)\times n(y))\cdot \varphi(x,y)\,dxd\sigma(y)= \int_{\Omega\times Y^*}f(x)\cdot  \varphi(x,y)\,dxdy\,.
\end{array}
$$
By density, this variational formulation holds for $\varphi\in \mathcal{V}$ with ${\rm div}_y(\varphi)=0\ \hbox{in }\Omega\times Y^*$, where
$$\mathcal{V}=
\Big\{\begin{array}{l}
 \varphi\in L^2(\Omega;H^1_{\rm per}(Y^*)^3)\ :\  \varphi\cdot n=0\ \hbox{ on }\Omega\times \partial F
 \end{array}\Big\}.$$
 Then we easily find that the function $\hat u$ satisfies the variational formulation 
\begin{equation}\label{form_var_hatu}\begin{array}{l}\displaystyle
 \int_{\Omega \times Y^*} \rot_y(\hat u(x,y))\cdot \rot_y(\varphi(x,y))dxdy
-2N^2\int_{\Omega \times Y^*}  \rot_{y}(\varphi(x,y))\cdot \hat w(x,y)dxdy
\\
\noame\displaystyle
\qquad
+2 \left({1\over \alpha}-N^2\right)\int_{\Omega\times \partial F}(\hat w(x,y)\times n(y))\cdot \varphi(x,y)\,dxd\sigma(y)= \int_{\Omega\times Y^*}f(x)\cdot  \varphi(x,y)\,dxdy\,,
\end{array}
\end{equation}
for every $\varphi\in \mathcal{V}^0(Y^*)$ where 
$$\mathcal{V}^0(Y^*)=
\left\{\begin{array}{l}
 \varphi\in \mathcal{V}\ :\  {\rm div}_y(\varphi)=0\ \hbox{in }\Omega\times Y^*,\\
 \noame
\displaystyle  {\rm div}_x\left(\int_{Y^*}\varphi\,dy\right)=0\ \hbox{in }\Omega,\quad \left(\int_{Y^*}\varphi\,dy\right)\cdot n=0\ \hbox{on }\partial\Omega
 \end{array}\right\}.$$

Next, we prove (\ref{system_1_3_hat}).  To do this, we take as test function in (\ref{var_form_v}) the function $\psi_\ep(x)= \eta(x) \Psi(x/\ep)$, where $\eta\in  \mathcal{D}(\Omega)$ and $\Psi \in H^1_{\rm per}(Y^*)^3$  with $\Psi\cdot n=0$ on $\partial F$, and we have
$$\begin{array}{l}\displaystyle
\ep^2R_c\int_{\Omega_\varepsilon}\rot(w_\ep)\cdot \rot(\psi_\ep)\,dx+\ep^2R_c\int_{\Omega_\varepsilon}{\rm div}(w_\ep)\cdot {\rm div}(\psi_\ep)\,dx+4N^2\int_{\Omega_\varepsilon}w_\ep\cdot\psi_\ep\,dx\\
\noame\displaystyle
\qquad-2N^2\int_{\Omega_\varepsilon}\rot(\psi_\ep)\cdot u_\ep\,dx + 2N^2(\beta-1)\int_{\partial F_\varepsilon}(u_\ep\times n_\ep)\cdot \psi_\ep\,d\sigma
= \int_{\Omega_\ep}g\cdot \psi_\ep\,dx\,.
\end{array}
$$
Similarly to what happens with  $\varphi_\ep$ and taking into account that 
\begin{equation}\label{property_psi_2}\rot(\psi_\ep)= \nabla\eta\times\Psi\left({\cdot\over \ep}\right)+ \ep^{-1}\eta\, \rot_y\left(\Psi\left({\cdot\over \ep}\right)\right), \quad {\rm div}(\psi_\ep)=\nabla\eta\cdot\Psi\left({\cdot\over \ep}\right)+ \ep^{-1}\eta\, {\rm div}_y\left(\Psi\left({\cdot\over \ep}\right)\right),\end{equation}
we have $\mathcal{T}_\ep(\psi_\ep)= \mathcal{T}_\ep(\eta) \Psi$ and 
\begin{equation}\label{property_psi}
\begin{array}{c}
\mathcal{T}_\ep(\psi_\ep)\to  \eta \Psi\quad\hbox{in }L^2(\Omega\times Y^*)^3,\\
\noame
\ep\mathcal{T}_\ep({\rm div}(\psi_\ep))
\to \eta\, {\rm div}_y(\Psi)\hbox{ in }L^2(\Omega\times Y^*)^3,\quad \ep\mathcal{T}_\ep(\rot(\psi_\ep))
\to \eta\, \rot_y(\Psi)\hbox{ in }L^2(\Omega\times Y^*)^3.
\end{array}
\end{equation}
By applying the unfolding to the variational formulation and using Propositions \ref{properties_UO}$_{2,3}$ and \ref{properties_BUO}$_{2,3,4}$, we get
$$\begin{array}{l}\displaystyle
\ep^2R_c\int_{\mathbb{R}^3\times Y^*}\mathcal{T}_\ep(\rot(w_\ep))\cdot \mathcal{T}_\ep(\rot(\psi_\ep))\,dxdy+\ep^2R_c\int_{\mathbb{R}^3\times Y^*}\mathcal{T}_\ep({\rm div}(w_\ep))\cdot \mathcal{T}_\ep({\rm div}(\psi_\ep))\,dxdy\\
\noame\displaystyle
\qquad+4N^2\int_{\mathbb{R}^3\times Y^*}\mathcal{T}_\ep(w_\ep)\cdot
\mathcal{T}_\ep(\psi_\ep)\,dxdy -2N^2\int_{\mathbb{R}^3\times Y^*}\mathcal{T}_\ep(\rot(\psi_\ep))\cdot \mathcal{T}_\ep(u_\ep)\,dxdy \\
\noame\displaystyle
\qquad+ 2N^2(\beta-1)\ep^{-1}\int_{\mathbb{R}^3\times \partial F}
\left(\mathcal{T}_\ep^b(u_\ep)\times n\right)\cdot \mathcal{T}_\ep(\psi_\ep)\,dxd\sigma(y)
= \int_{\mathbb{R}^3\times Y^*}\mathcal{T}_\ep(g)\cdot \mathcal{T}_\ep(\psi_\ep)\,dxdy\,.
\end{array}
$$
Next, we pass to the limit in every terms of the previous variational formulation:
\begin{itemize}
\item First to third terms. From convergences (\ref{conv_hat_w_1}), (\ref{conv_hat_w_1_rot}) and (\ref{property_psi}), we have 
$$\begin{array}{l}
\displaystyle \ep^2\int_{\mathbb{R}^3\times Y^*}\mathcal{T}_\ep(\rot(w_\ep))\cdot \mathcal{T}_\ep(\rot(\psi_\ep))\,dxdy\displaystyle 
=\int_{\mathbb{R}^3\times Y^*}\ep \mathcal{T}_\ep(\rot(w_\ep))\cdot \ep \mathcal{T}_\ep(\rot(\psi_\ep))\,dxdy\\
\noame
\displaystyle
\qquad \to \int_{\Omega\times Y^*}\rot_y(\hat w(x,y))\cdot  (\eta(x)\rot_y(\Psi(y)))\,dxdy=\int_{\Omega\times Y^*}\rot_y(\hat w(x,y))\cdot  \rot_y(\eta(x)\Psi(y))\,dxdy\,,
\end{array}$$
$$\begin{array}{l}
\displaystyle \ep^2\int_{\mathbb{R}^3\times Y^*}\mathcal{T}_\ep({\rm div}(w_\ep))\cdot \mathcal{T}_\ep({\rm div}(\psi_\ep))\,dxdy=
\int_{\mathbb{R}^3\times Y^*}\ep \mathcal{T}_\ep({\rm div}(w_\ep))\cdot \ep \mathcal{T}_\ep({\rm div}(\psi_\ep))\,dxdy\\
\noame
\qquad \displaystyle \to 
\int_{\Omega\times Y^*} {\rm div}_y(\hat w(x,y))\cdot (\eta(x){\rm div}_y(\Psi(y)))\,dxdy=\int_{\Omega\times Y^*} {\rm div}_y(\hat w(x,y))\cdot {\rm div}_y(\eta(x) \Psi(x,y))\,dxdy,
\\
\\
\displaystyle \int_{\mathbb{R}^3\times Y^*}\mathcal{T}_\ep(w_\ep)\cdot
\mathcal{T}_\ep(\psi_\ep)\,dxdy\to \int_{\Omega\times Y^*}\hat w(x,y)\cdot (\eta(x)\Psi(y))\,dxdy.
\end{array}$$
\item Fourth term. From (\ref{conv_hat_u_1}) and (\ref{property_psi}), we get
$$\begin{array}{rl}
\displaystyle\int_{\mathbb{R}^3\times Y^*}\!\!\mathcal{T}_\ep(\rot(\psi_\ep))\cdot \mathcal{T}_\ep(u_\ep)\,dxdy=&\displaystyle\!\!\!
\int_{\mathbb{R}^3\times Y^*}\!\!\ep\mathcal{T}_\ep(\rot(\psi_\ep))\cdot \ep^{-1}\mathcal{T}_\ep(u_\ep)\,dxdy
\\
\noame
\displaystyle
\to &\displaystyle\!\!\!\!\int_{\Omega\times Y^*}\!\!\eta(x)\rot_y(\Psi(y))\cdot \hat u(x,y)\,dxdy=\int_{\Omega\times Y^*}\!\!\rot_y(\eta(x)\Psi(y))\cdot \hat u(x,y)\,dxdy.
\end{array}$$
\item Fifth term. From convergence (\ref{conv_hat_u_1_b}) and (\ref{property_psi}),
$$\ep^{-1}\int_{\mathbb{R}^3\times \partial F}
\left(\mathcal{T}_\ep^b(u_\ep)\times n(y))\right)\cdot \mathcal{T}_\ep(\psi_\ep)\,dxd\sigma(y)\to 
\int_{\Omega\times \partial F}
\left(\hat u(x,y)\times n(y)\right)\cdot (\eta(x) \Psi(y))\,dxd\sigma(y).$$
\item Sixth term. From Proposition \ref{properties_UO}$_7$ and convergence (\ref{property_psi}), 
$$\int_{\mathbb{R}^3\times Y^*}\mathcal{T}_\ep(g)\cdot \mathcal{T}_\ep(\psi_\ep)\,dxdy\to \int_{\mathbb{R}^3\times Y^*}g(x)\cdot (\eta(x)\Psi(y))\,dxdy.$$
\end{itemize}

From the previous convergences and noting $\psi(x,y)=\eta(x)\Psi(y)$, we obtain 
\begin{equation}\label{form_var_hatw}\begin{array}{l}\displaystyle
R_c\int_{\Omega\times Y^*} \rot_y(\hat w(x,y))\cdot  \rot_y(\psi(x,y))\,dxdy+R_c\int_{\Omega\times Y^*} {\rm div}_y(\hat w(x,y))\cdot {\rm div}_y( \psi(x,y))\,dxdy\\
\noame\displaystyle
\qquad+4N^2\int_{\Omega\times Y^*}\hat w(x,y)\cdot
 \psi(x,y)\,dxdy
-2N^2\int_{\Omega\times Y^*} \rot_y(\psi(x,y))\cdot \hat u(x,y)\,dxdy \\
\noame\displaystyle
\qquad+ 2N^2(\beta-1)\int_{\Omega\times \partial F}
\left(\hat u(x,y)\times n(y)\right)\cdot  \psi(x,y)\,dxd\sigma(y)
= \int_{\Omega\times Y^*}g(x)\cdot \psi(x,y)\,dxdy\,,
\end{array}
\end{equation}
which, by density,  holds for $\hat \psi\in \mathcal{V}$.

Finally, under condition (\ref{condition_existence}), we can prove that the variational formulation given by (\ref{form_var_hatu}) and  (\ref{form_var_hatw})  admits a unique solution $(\hat u,\hat w)$ in $\mathcal{V}^0(Y^*)\times \mathcal{V}$ (see proof of Lemma \ref{lema_cell_problems} above), and so we conclude that the whole sequence $(\ep^{-1}\mathcal{T}_\ep(u_\ep), \mathcal{T}_\ep(w_\ep))$ converges to this solution. Reasoning as in \cite{Allaire0}, it follows that there exists $q(x)\in L^2(Y^*)/\mathbb{R}$ and $\hat q(x,y)\in L^2(\Omega;L^2_{{\rm per}}(Y^*)/\mathbb{R})$  such that  the variational formulation given by (\ref{form_var_hatu}) and (\ref{form_var_hatw}) is equivalent to system (\ref{system_1_1_hat})-(\ref{BC_holes_4_hat}).  It remains to prove that $q$ coincides with pressure $p$. This can be easily done by proceeding as above by considering a test function which is divergence-free only in $y$, and identifying limits.

{\it Step 2}. In this step, we eliminate the microscopic variable $y$ in the homogenized
problem (\ref{system_1_1_hat})-(\ref{BC_holes_4_hat}) and then, we obtain a Darcy equation for the pressure  $p$.
To do that, for every $k=1, 2$, $i=1, 2, 3$, we consider the micropolar local problem (\ref{local_problems}), which has  a unique solution  $(u^{i,k}, w^{i,k}, \pi^{i,k})\in V^0_Y\times V_Y\times L^2_0(Y^*)$ (see Lemma  \ref{lema_cell_problems} above). Then, by the  following indentification
$$\begin{array}{l}
\displaystyle\hat u(x,y)=\sum_{i=1}^3\left(\left(f_i(x)-\partial_{x_i}p(x)\right)u^{i,1}(y)+g_i(x)u^{i,2}(y)\right),
\end{array}$$
$$\begin{array}{l}
\displaystyle\hat w(x,y)=\sum_{i=1}^3\left(\left(f_i(x)-\partial_{x_i}p(x)\right)w^{i,1}(y)+g_i(x)w^{i,2}(y)\right),
\end{array}$$
$$\begin{array}{l}
\displaystyle\hat q(x,y)=\sum_{i=1}^3\left(\left(f_i(x)-\partial_{x_i}p(x)\right)\pi^{i,1}(y)+g_i(x)\pi^{i,2}(y)\right),
\end{array}$$
and thanks to identities (\ref{relation_u_hatu}) and (\ref{relation_w_hatu}), we deduce that $u$ and $w$ have  the expressions given in  (\ref{Darcy_law_vel_mic}). 

Next, the divergence condition with respect to the variable $x$ given in (\ref{divxhatu}) together with the expression of $u$ gives (\ref{Darcy_law}).

Finally,  positive definiteness of $K^{(1)}$ follows from the fact that $K^{(1)}_{ij}=\mathcal{A}_Y(u^{i,1},w^{i,1};u^{j,1},w^{j,1})$ and the coercivity of the bilinear form $\mathcal{A}_Y$ defined by (\ref{Def_A_cell}).  Since (\ref{Darcy_law}) is an elliptic equation with   $K^{(1)}f+K^{(2)}g\in L^2(\Omega)$, then it has a unique solution $p\in H^1(\Omega)\cap L^2_0(\Omega)$, and $u, w\in L^2(\Omega)^3$ given by (\ref{Darcy_law_vel_mic}) are also unique. By the uniqueness of the limits, the whole sequences converge. 
\end{proof}

%\section*{Conflict of interest}
%This work does not have any conflicts of interest.
%
%\section*{Acknowledgement}
%There are no funders to report for this submission.

\end{document}